\documentclass[10pt]{article}  
\usepackage[english]{babel}
\usepackage[toc,page]{appendix}
\usepackage{bbm}
\usepackage{fullpage}
\usepackage[top=1in, bottom=1in, left=1in, right=1in]{geometry}
\usepackage{etaremune}
\usepackage{enumitem}

\usepackage{amssymb,amsmath,amsthm}

\usepackage{xcolor}
\usepackage{hyperref}

\usepackage[showonlyrefs]{mathtools}

\DeclareSymbolFont{rsfs}{U}{rsfs}{m}{n}
\DeclareSymbolFontAlphabet{\mathscrsfs}{rsfs}
\usepackage{mathrsfs}

\usepackage{url}

\hypersetup{
  colorlinks   = true, 
  urlcolor     = blue, 
  linkcolor    = blue, 
  citecolor   = red 
}

\usepackage[labelformat=empty]{caption}

\newtheorem{theorem}{Theorem}[section]
\newtheorem{lemma}[theorem]{Lemma}

\newtheorem{proposition}[theorem]{Proposition}
\newtheorem{corollary}[theorem]{Corollary}

\theoremstyle{definition}

\newtheorem{remark}[theorem]{Remark}

\numberwithin{equation}{section}

\usepackage{times}

\newcommand{\bea}{\begin{eqnarray}}
\newcommand{\eea}{\end{eqnarray}}
\newcommand{\<}{\langle}
\renewcommand{\>}{\rangle}

\newcommand{\wt}{\widetilde}
\newcommand{\op}{\text{op}}

\def\generic{{\rm generic}}

\def\eigen{{\rm eigen}}
\def\bdd{{\rm bounded}}

\def\iid{\text{i.i.d.~}}

\def\eps{{\varepsilon}}

\def\supp{{\rm supp}}

\def\hq{{\widehat{q}}}

\def\bx{{\boldsymbol{x}}}

\def\cT{{\mathcal T}}
\def\cC{{\mathcal C}}

\def\blambda{{\boldsymbol \lambda}}

\def\op{{\rm op}}

\def\bsig{{\boldsymbol {\sigma}}}

\def\by{{\boldsymbol y}}

\def\bv{{\boldsymbol{v}}}

\def\bx{{\boldsymbol{x}}}

\def\bA{\boldsymbol{A}}
\def\bB{\boldsymbol{B}}

\def\bL{\boldsymbol{L}}

\def\de{{\rm d}}

\def\<{\langle}
\def\>{\rangle}
\def\Tr{{\sf Tr}}

\def\rank{{\rm rank}}

\def\cM{{\cal M}}
\def\cN{{\cal N}}

\def\by{{\boldsymbol{y}}}

\def\bw{{\boldsymbol{w}}}

\def\blambda{{\boldsymbol{\lambda}}}

\def\b0{{\boldsymbol{0}}}

\def\rd{{\mathrm {rad}}}

\def\bG{{\boldsymbol G}}

\DeclareMathOperator*{\plim}{p-lim}

\def\cS{{\mathcal S}}

\def\cK{{\mathcal K}}

\renewcommand{\b}{\mathbf{b}}

\def\fr{\frac}
\def\lt{\left}
\def\rt{\right}

\def\la{\langle}
\def\ra{\rangle}

\def\eps{\varepsilon}

\def\bbE{{\mathbb{E}}}

\def\bbN{{\mathbb{N}}}
\def\bbP{{\mathbb{P}}}
\def\bbR{{\mathbb{R}}}

\def\bbW{{\mathbb{W}}}

\def\cK{{\mathcal{K}}}
\def\cN{{\mathcal{N}}}
\def\cP{{\mathcal{P}}}
\def\cQ{{\mathcal{Q}}}

\def\sH{{\mathscr{H}}}

\def\GS{{\mathrm{GS}}}

\DeclareMathOperator*{\E}{\bbE}

\DeclareMathOperator*{\argmin}{\rm arg\,min}

\def\sph{\mathrm{sp}}

\newcommand{\Gp}[1]{\mathbf{G}^{(#1)}}

\newcommand{\norm}[1]{{\lt\|#1\rt\|}}
\newcommand{\tnorm}[1]{{\|#1\|}}

\newcommand{\rev}[1]{{#1}}
\newcommand{\revv}[1]{{#1}}

\author{Mark Sellke\thanks{Harvard University, Department of Statistics. Cambridge, MA, USA. Email: \texttt{msellke@fas.harvard.edu}}}

\title{On Marginal Stability in Low Temperature Spherical Spin Glasses}

\date{}

\begin{document}

\maketitle

\begin{abstract}
\noindent
We show marginal stability of near-ground states in spherical spin glasses is equivalent to full replica symmetry breaking at zero temperature near overlap $1$.
This connection has long been implicit in the physics literature, which also links marginal stability to the performance of efficient algorithms.
For even models, we rule out Hessian eigenvalues above the upper edge of the bulk spectrum, and obtain geometric consequences for low temperature Gibbs measures in the case that marginal stability is absent.
Our proofs rely on interpolation bounds for vector spin glass models.
For generic models, we give another more conceptual argument that full RSB near overlap $1$ implies marginal stability at low temperature.
\end{abstract}

\section{Introduction}

Spherical spin glass Hamiltonians are random disordered smooth functions in very high dimension.
Their landscapes are understood to be rich and complicated, with exponentially many local maxima at a range of energy levels.
We will be interested in the qualitative behavior around their near-global maxima.
Does local uniform concavity hold near the extreme values, or could the Hessian be ill-conditioned?
The former would imply that low temperature Gibbs measures are supported on isolated wells, and that the corresponding stationary dynamics remain trapped within them.
By contrast an ill-conditioned \emph{marginally stable} Hessian would allow for the possibility of a connected ``manifold'' of near maxima.
In the physics literature, marginal stability at low-temperature is widely believed to be equivalent to \emph{full replica symmetry breaking} (full RSB) near overlap $1$ at zero temperature, which is a property of the order parameter in the Parisi formula.
We prove a strong form of this equivalence for spherical spin glasses with even interactions, and derive consequences for Langevin dynamics and disorder chaos whenever full RSB is absent.

As our results depend on properties of the minimizer in Parisi's variational formula, we begin by recalling this formula.
For each $p\geq 1$, let $\Gp{p} \in \lt(\bbR^N\rt)^{\otimes p}$ be an independent $p$-tensor with \iid standard Gaussian entries. 
Fixing an infinite sequence $(\gamma_p)_{p\geq 1}$ of non-negative reals, the mixed $p$-spin Hamiltonian $H_N$ is 
\begin{align}
    \label{eq:def-hamiltonian}
    H_N(\bsig) &= \sum_{p\geq 1} \fr{\gamma_p}{N^{(p-1)/2}} \la \Gp{p}, \bsig^{\otimes p} \ra.
\end{align}
The coefficients $\gamma_p$ are encoded in the \emph{mixture function} $\xi(x) = \sum_{p\geq 1} \gamma_p^2 x^p$, which we assume is not linear and has radius of convergence strictly larger than $1$. 
We view $H_N$ as a function on the spherical domain $\cS_N \equiv \big\{\bsig \in \bbR^N : \sum_{i=1}^N \bsig_i^2 = N\big\}$; equivalently, it is the centered Gaussian process on $\cS_N$ with covariance
\begin{equation}
\label{eq:covariance-formula}
    \E H_N(\bsig^1) H_N(\bsig^2) = N \xi(R(\bsig^1,\bsig^2))
    \equiv 
    N \xi(\la \bsig^1,\bsig^2\ra/N)
    .
\end{equation}
Here $R(\bsig^1,\bsig^2)=\la \bsig^1,\bsig^2\ra/N\in [-1,1]$ is known as the \emph{overlap} between $\bsig^1,\bsig^2$.

We will study $H_N$ near its extreme values, where the ground state energy
\begin{equation}
\begin{aligned}
    \GS_N=\max_{\bsig\in \cS_N} H_N(\bsig)/N.
\end{aligned}
\end{equation}
is approximately achieved.
The limiting value of $\GS_N$ is given by Parisi's formula.
Following \cite{chen2017parisi} (see also \cite{crisanti1992spherical,talagrand2006free,chen2013aizenman}), let $\cN$ denote the set of non-decreasing, right-continuous functions $\zeta:[0,1)\to \bbR_{\geq 0}$ equipped with the vague topology, and define
\begin{equation}
\label{eq:def-cK}
    \cK=\lt\{
    (\zeta,L)\in \cN\times(0,\infty)~:~L>\int_0^1 \zeta(s)~\de s
    \rt\}
    .
\end{equation}
For $(\zeta,L)\in\cK$, define $\hat\zeta(q)=L-\int_0^q \zeta(s)~\de s$ and:
\[
    \cQ(\zeta,L)
    =
    \frac{1}{2}
    \lt(
    \xi'(0)L + \int_0^1 \xi''(q)\hat\zeta(q)~\de q + \int_0^1 \frac{\de q}{\hat\zeta(q)}
    \rt)
    =
    \frac{1}{2}
    \lt(
    \xi'(1)L
    \rev{-}
    \int_0^1 \xi''(q)\Big(\int_0^q \zeta(s)\de s\Big)~\de q + \int_0^1 \frac{\de q}{\hat\zeta(q)}
    \rt)
    .
\]

\begin{proposition}[{\cite{chen2017parisi}}]
\label{prop:parisi-zero-temp}
    The in-probability limit of the ground state energy is:
    \begin{equation}
    \label{eq:cs-functional-0temp-inf}
    \GS(\xi)
    \equiv
    \plim_{N\to\infty} \GS_N
    =
    \inf_{(\zeta,L)\in\cK} \cQ(\zeta,L)
    .
    \end{equation}
    Further, there exists a unique minimizing pair $(\zeta,L)\in\cK$.
\end{proposition}

The minimizing $(\zeta,L)$ is characterized by local stationarity conditions, which are reviewed in Subsection~\ref{subsec:characterization-of-minimizer}.

\subsection{Main Results}

As previously mentioned, we study the local behavior of $H_N$ around its extreme values.
To start, Proposition~\ref{prop:main} shows the bulk Hessian spectrum at the ground state is directly described by the minimizer $(\zeta,L)\in\cK$ in the zero-temperature Parisi formula \eqref{eq:cs-functional-0temp-inf}, in particular the value $\hat\zeta(1)=L-\int_0^1 \zeta(s)\de s$.
Namely the rescaled radial derivative and bulk spectral edges of the Hessian are asymptotically given by the formulas:
\begin{align}
\label{eq:r-def}
    r(\xi)
    &=
    \hat\zeta(1)\xi''(1)+\hat\zeta(1)^{-1},
    \\
\label{eq:lambda-def}
    \lambda_{\pm}(\xi)
    &=
    \pm 2\sqrt{\xi''(1)}
    -
    r(\xi)
    =
    -
    \hat\zeta(1)
    \lt(
    \sqrt{\xi''(1)}
    \mp
    \hat\zeta(1)^{-1}
    \rt)^2
    \leq 0.
\end{align}
The same description extends, up to error $o_{\delta\to 0}(1)$, to all \emph{$\delta$-approximate ground states}, i.e. points in the set 
\begin{equation}
\label{eq:approx-ground-states}
    A_{\delta}
    =
    \big\{\bsig\in\cS_N~:~H_N(\bsig)/N\geq \GS(\xi)-\delta\big\}
\end{equation}
Below and throughout, we write $\blambda_{k}(\cdot)$ for the $k$-th largest eigenvalue of a symmetric matrix.
$\nabla_{\sph}^2 H_N(\cdot)$ denotes the Riemannian Hessian on $\cS_N$, an $(N-1)\times (N-1)$ matrix defined in Subsection~\ref{subsec:prelims}.
We often say an event depending on $H_N$ and other parameters (e.g.\ $(\xi,\eps,\delta)$) holds with probability (at least) $1-e^{-cN}$.
Here $c$ is always sufficiently small depending on the other parameters, while $N$ is sufficiently large depending on everything else including $c$.

\begin{proposition}
    \label{prop:main}
    For any $\xi$ and $\eps>0$, there exists $\delta>0$ such that the following holds with probability $1-e^{-cN}$ for large $N$.
    For all $\delta$-approximate ground states $\bsig\in A_{\delta}$, the \rev{Euclidean} gradient of $H_N$ \rev{(in $\bbR^N$)} satisfies
    \begin{equation}
        \label{eq:radial-deriv-formula}
        \|\nabla H_N(\bsig) 
        -
        r(\xi)\,\bsig\|/\sqrt{N}
        \leq 
        \eps
    \end{equation}
    while the top and bottom of the bulk Hessian spectrum satisfy
    \begin{equation}
    \begin{aligned}
    \label{eq:bulk-edge-formula}
        \Big|\blambda_{\lfloor\delta N\rfloor}\big(\nabla_{\sph}^2 H_N(\bsig)\big)
        -\lambda_+(\xi)
        \Big|
        &\leq 
        \eps,
        \\
        \Big|\blambda_{N-\lfloor\delta N\rfloor}\big(\nabla_{\sph}^2 H_N(\bsig)\big)
        -\lambda_-(\xi)
        \Big|
        &\leq
        \eps
        .
    \end{aligned}
    \end{equation}
\end{proposition}

This result is stated as a proposition because a more abstract formula for $\nabla H_N(\bsig)$, expressed as a derivative of the function $\xi\mapsto \GS(\xi)$, is essentially in \cite[Remark 2]{chen2017parisi} (or see \cite[Corollary 7]{subag2018free}).
Although we provide a detailed proof for completeness, the only novelty is an integration by parts using the stationarity conditions for $(\zeta,L)$, which yields the more tractable \eqref{eq:r-def}.
Regarding the bulk spectrum \eqref{eq:bulk-edge-formula}, known concentration estimates imply that uniformly for all $\bsig\in\cS_N$, the spectral measure of $\nabla_{\sph}^2 H_N(\bsig)$ is approximately a semi-circle density with radius $2\sqrt{\xi''(1)}$ and shifted by the rescaled radial derivative $\partial_{\rd}H_N(\bsig)=R(\bsig,\nabla H_N(\bsig))$; see the proof of \cite[\rev{Lemma 1.3}]{subag2018following} (which is rephrased below as Lemma~\ref{lem:spectrum-approx}).
Thus \eqref{eq:bulk-edge-formula} follows routinely from \eqref{eq:radial-deriv-formula}.

Note that $\lambda_+$ defined above is never positive, and equals $0$ if and only if
\begin{equation}
\label{eq:full-RSB-behavior}
\hat\zeta(1)=\xi''(1)^{-1/2}.
\end{equation}
We say $\xi$ exhibits \emph{full RSB endpoint behavior} when this condition holds.
As explained in Corollary~\ref{cor:full-RSB-justification}, this condition is implied by (and is presumably generically equivalent to) the usual definition of full RSB that $\zeta$ is strictly increasing in a neighborhood of $1$. 
\revv{
We occasionally refer to the latter as \emph{strictly full RSB endpoint behavior}.
}

When \eqref{eq:full-RSB-behavior} holds, Proposition~\ref{prop:main} states that the Hessian at any approximate ground state has many near-zero eigenvalues, i.e\ is marginally stable.
Conversely when $\lambda_+<0$, it is natural to hope that $H_N$ is locally uniformly concave at near-extrema.
However it is not \rev{at all obvious from existing results that the Hessian does not have a constant number of eigenvalues above the upper edge of the bulk spectrum}. 
Our next result rules out such upward outliers when $\xi$ contains only even degree terms.\footnote{
Theorem~\ref{thm:no-upward-outliers-even-xi} and Corollary~\ref{cor:marginal-stability} extend with no changes when $\gamma_1>0$, but do require $\gamma_p=0$ for odd $p\geq 3$.
The same holds for Corollaries~\ref{cor:deep-level-sets} and \ref{cor:langevin}. 
However deducing slow mixing from Corollary~\ref{cor:langevin} requires modification if $\gamma_1>0$, as does Corollary~\ref{cor:disorder-chaos}. 
Indeed for the former, one must exclude the topologically trivial phase where low-temperature Langevin dynamics does mix rapidly \cite[Theorem 1.8]{huang2023strong}.
}
As discussed further in Subsection~\ref{subsec:proof-ideas}, the idea is to connect the largest Hessian eigenvalue to the ground state energy of a vector spin glass, which can be controlled using a precise interpolation bound.

\begin{theorem}
\label{thm:no-upward-outliers-even-xi}
    Suppose $\xi$ is even, i.e. $\gamma_p=0$ for all odd $p$. 
    Then for any $\eps$ there is $\delta$ such that the following holds with probability $1-e^{-cN}$ for large $N$.
    For all $\bsig\in A_{\delta}$,
    \begin{equation}
        \label{eq:top-edge-formula}
        \Big|\blambda_1\big(\nabla_{\sph}^2 H_N(\bsig)\big)
        -
        \lambda_+(\xi)
        \Big|
        \leq 
        \eps.
    \end{equation}
\end{theorem}

The following corollary is essentially immediate.
(For the easy upper bound $\blambda_1\leq \eps/2$ in \eqref{eq:marginal-stability} without assuming $\xi$ is even, see Corollary~\ref{cor:near-ground-state}.)

\begin{corollary}
    \label{cor:marginal-stability}
    If $\xi$ exhibits full RSB endpoint behavior \eqref{eq:full-RSB-behavior}, then all approximate ground states are marginally stable.
    Namely for any $\eps>0$, if $\delta=\delta(\xi,\eps)$ is small enough, then with probability $1-e^{-cN}$ all $\bsig\in A_{\delta}$ satisfy:
    \begin{equation}
    \label{eq:marginal-stability}
    |\blambda_{\lfloor\delta N\rfloor}\big(\nabla_{\sph}^2 H_N(\bsig)\big)|
    +
     |\blambda_{1}\big(\nabla_{\sph}^2 H_N(\bsig)\big)|
    \leq 
    \eps.
    \end{equation}
    Conversely, suppose $\xi$ does not exhibit full RSB endpoint behavior and is even, and let $\bsig\in A_{\delta}$ for small $\delta\leq\delta_*(\xi)$.
    Then $H_N$ is locally uniformly concave near $\bsig$ in that with probability $1-e^{-cN}$:
    \begin{equation}
    \label{eq:local-concavity}
    \blambda_1\big(\nabla_{\sph}^2 H_N(\bsig)\big)
    \leq 
    -1/C(\xi)
    <
    0.
    \end{equation}
\end{corollary}

We note that these results apply to Gibbs samples at temperature tending to $0$ slowly with $N$, simply because they are approximate ground states (see e.g.\ Lemma~\ref{lem:small-values-low-prob}). 
They also have consequences at positive temperatures not tending to $0$.
Here as usual, the Gibbs measure $\mu_{\beta}=\mu_{\beta,H_N}$ at inverse temperature $\beta$ is defined by
\[
\de\mu_{\beta}(\bsig)=e^{\beta H_N(\bsig)}\de\mu_0(\bsig)/Z_{N,\beta}.
\]
$Z_{N,\beta}=\int_{\cS_N} e^{\beta H_N(\bsig)} \de \mu_0(\bsig)$ is the partition function relative to uniform measure $\mu_0$ on $\cS_N$.
Under genericity conditions on $\xi$, $\mu_{\beta}$ is arranged into an ultrametric tree of \emph{ancestor states}, which are certain points in the interior of $\cS_N$ that are approximate ground states for their respective radius $t\sqrt{N}$; see \cite{panchenko2013parisi,jagannath2017approximate,chatterjee2021average,subag2018free}. 
Then applying the above results to $H_N|_{t\cS_N}$, which due to the formula \eqref{eq:covariance-formula} amounts to studying $\xi_t(q)=\xi(t^2 q)$, describes the local behavior of $H_N$ near ancestor states. 
The corresponding zero-temperature order parameter $(\zeta_t,L_t)$ can be directly read off from the positive-temperature analog of the original model, and will exhibit full RSB endpoint behavior if and only if the positive-temperature Parisi order parameter for $\xi$ exhibits analogous full RSB behavior at overlap $t^2$ (see \cite[Proposition 11]{subag2018free}).

Further, the local uniform concavity in \eqref{eq:local-concavity} easily implies a qualitative description of the deep level sets $A_{\delta}$.
As stated in Corollary~\ref{cor:deep-level-sets} below, all connected components $A_{\delta}$ consist of small separated clusters; this is a more geometric formulation of non-marginal stability.

\begin{corollary}
\label{cor:deep-level-sets}
    Suppose $\xi$ does not exhibit full RSB endpoint behavior and is even.
    For some $C(\xi)>0$ and all small enough $\delta\in (0,\delta_*(\xi))$, let $A_{\delta}$ be as in \eqref{eq:approx-ground-states}.
    Then with probability $1-e^{-cN}$, every connected component of $A_{\delta}$ contains exactly $1$ critical point of $H_N$, which is a local maximum, and has diameter at most $\sqrt{N\delta}/C(\xi)$.
    Further, distinct components of $A_{\delta}$ are separated by $\sqrt{N}/C(\xi)$.
\end{corollary}

\revv{
A similar result is deduced in \cite[Proposition 9.1]{arous2020geometry} for a special class of $\xi$, for which a version of \eqref{eq:local-concavity} follows from Kac--Rice asymptotics.
We therefore give only a proof outline (just below), following their approach.
As in \cite{arous2020geometry} our outline uses Hessian control only at \emph{exact} critical points, so it may be simplifiable.

\begin{proof}[Proof Outline of Corollary~\ref{cor:deep-level-sets}]
Consider spherical gradient-ascent flow $(\bsig_t)_{t\geq 0}$ starting from any $\bsig\in A_{\delta}$.
This flow must have an accumulation point on the sphere by compactness, and monotonicity of $H_N(\bsig_t)$ implies accumulation occurs at some $\bsig_*\in A_{\delta}$.
The same monotonicity also implies that $\bsig_*$ is a critical point for $H_N$. 
Then $\blambda_1(\nabla_{\sph}^2 H_N(\bsig_*))\leq -C(\xi)<0$ by \eqref{eq:local-concavity}, so $\bsig_*$ is a local maximum.

We now show the remaining geometric claims.
Choose $C'>0$ small compared to $C$ and $C''$ large compared to $C$, all depending only on $\xi$ (and in particular not on $\delta$).
Standard smoothness estimates on $H_N$ (see Proposition~\ref{prop:gradients-bounded}) then imply that the intersection $\cC(\bsig_*)$ of $A_{\delta}$ with the ball $B_{C'\sqrt{N}}(\bsig_*)$ contains $B_{C'\sqrt{\delta N}}(\bsig_*)$ and is contained within $B_{C''\sqrt{\delta N}}(\bsig_*)$.
Further, $\cC(\bsig_*)$ is connected and thus is a connected component of $A_{\delta}$ (indeed $\cC(\bsig_*)$ is star-shaped relative to the center point $\bsig_*$ since $\nabla_{\sph}^2 H_N$ is nearly locally constant).
This yields the diameter claim (modulo adjusting $C$).
Since $\bsig_t$ is a continuous path staying within $A_{\delta}$ it follows that $\bsig\in \cC(\bsig_*)$.
Since $\delta$ was small enough depending on $\xi$, we have $C/2\geq C''\sqrt{\delta}$.
Hence the $C\sqrt{N}/2$-neighborhood of $\cC(\bsig_*)$ contains no points in $A_{\delta}\backslash \cC(\bsig_*)$, which yields the separation claim (modulo adjusting $C$).
\end{proof}
}


We pause here to emphasize two points.
First, Corollary~\ref{cor:deep-level-sets} relies crucially on the absence of upward Hessian outliers stated in \eqref{eq:top-edge-formula}, as even a single near-zero upward outlier eigenvalue could result in $A_{\delta}$ having large components. 
Second, the full RSB assumption is essentially necessary: when $\lambda_+=0$, Corollary~\ref{cor:deep-level-sets} evidently does not hold at any $\bsig\in A_{\delta}$.
Indeed, it is predicted in e.g.\ \cite[Section 9]{cugliandolo1994out} that under full RSB behavior,  $\delta=O(1/N)$ may suffice for $A_{\delta}$ to have components of macroscopic diameter.
Rigorously, it follows (e.g.\ by generic perturbations as in \cite[Chapter 3.7]{panchenko2013sherrington}) that when $\xi$ is full RSB near overlap $1$ in the stronger sense that $\zeta$ from \eqref{eq:cs-functional-0temp-inf} is strictly increasing on $[q,1]$, then for any $\delta>0$ and sufficiently large $N$, with high probability there exist $\bsig,\bsig'\in A_{\delta}$ with $R(\bsig,\bsig')=q$.


Corollary~\ref{cor:deep-level-sets} yields further consequences for low temperature Gibbs measures in the absence of marginal stability, \revv{which stem from similarities with the shattered phase studied in \cite{alaoui2023shattering}}.
\revv{
In fact a substantial portion of Corollary~\ref{cor:deep-level-sets} holds slightly more generally, in the absence of \emph{strictly} full RSB endpoint behavior (defined below \eqref{eq:full-RSB-behavior}).
As explained in Remark~\ref{rem:weaker-condition} at the end of this section, this also suffices for the two corollaries below.}

Using the separation of distinct components, and that low temperature Gibbs samples are approximate ground states, one obtains a plateau property for the autocorrelation function of stationary Langevin dynamics. 
Given an initialization $\bx_0\in\cS_N$, the spherical Langevin dynamics are given by 
\begin{equation}
\label{eq:langevin}
\de \bx_t=\lt(\beta \nabla_{\sph} H_N(\bx_t)-\frac{(N-1) \bx_t}{N}\rt)\de t + P_{\bx_t}^{\perp} \sqrt{2}\,\de \bB_t
\end{equation}
where $P_{\bx_t}^{\perp}=I_N-\bx_t^{\otimes 2}/N$ is a rank \rev{$N-1$} projection and $\bB_t$ is standard $N$-dimensional Browian motion.
It is well known that these dynamics remain on $\cS_N$ almost surely, with unique stationary distribution $\mu_{\beta}$.
We omit the proof of the corollary below, which is exactly the same as \cite[Corollary \rev{2.9}]{alaoui2023shattering}.
A byproduct is exponentially slow mixing of Langevin dynamics under the same conditions (since $\mu_{\beta}$ is origin-symmetric when $\xi$ is even).

\begin{corollary}
\label{cor:langevin}
    Suppose $\xi$ does not exhibit full RSB endpoint behavior and is even, and fix $\eps>0$.
    For some $C(\xi)>0$ and all large enough $\beta\geq \beta(\xi,\eps)$, let $\bx_0\sim\mu_{\beta}$ be a Gibbs sample, and $\bx_t$ be the trajectory of spherical Langevin dynamics \eqref{eq:langevin}. 
    Then with probability $1-e^{-cN}$, one has that
    \[
    \inf_{0\leq t\leq e^{cN}}
    R(\bx_0,\bx_t)\geq 1-\eps.
    \]
\end{corollary}

The last assertion of Corollary~\ref{cor:deep-level-sets} also allows us to determine the scale for the onset of transport disorder chaos, i.e.\ the sensitivity of $\mu_{\beta}$ to perturbations of $H_N$. 
Here we also assume $\xi$ is ``even generic'', with $\sum_{p\geq 1} \frac{1_{\gamma_p\neq 0}}{p}=\infty$.
With $\tilde H_N$ an independent copy of $H_N=H_{N,0}$, set
\[
H_{N,t}
=
\sqrt{1-t} H_N + \sqrt{t}\tilde H_N,\quad 
\forall\, t\in [0,1].
\]
Let $\mu_{\beta,t}=\mu_{\beta,H_{N,t}}$ be the corresponding Gibbs measure, and $\bbW_{2,N}(\mu,\mu')$ for the rescaled Wasserstein-$2$ distance
\[
W_{2,N}(\mu,\mu')^2
=
\inf_{\pi\in\Pi(\mu,\mu')}
\bbE^{\pi}[\|\bx-\by\|_2^2/N].
\]
Here the infimum is over couplings $(\bx,\by)\sim\pi$ with marginals $\bx\sim\mu$ and $\by\sim\mu'$.

\begin{corollary}
\label{cor:disorder-chaos}
    Suppose $\xi$ is even generic and does not exhibit full RSB endpoint behavior.
    Let $\beta\geq \beta(\xi)$ be sufficiently large.
    If $(\eps_N)_{N\geq 1}$ is a deterministic sequence with $\liminf_{N\to\infty} N\eps_N>0$, then
    \begin{equation}
    \label{eq:chaos}
    \liminf_{N\to\infty}
    \bbE[
    \bbW_{2,N}(\mu_{\beta,0},\mu_{\beta,\eps_N})
    ]
    >0.
    \end{equation}
\end{corollary}

The same $1/N$ scaling of $\eps_N$ was identified for pure spherical spin glasses at very low temperature \cite{subag2017geometry}, and established later within the shattered phase by \cite{alaoui2023shattering}.
Indeed as we explain in Subsection~\ref{subsec:chaos}, Corollary~\ref{cor:disorder-chaos} follows by the arguments in \cite[Section 5.1]{alaoui2023shattering}.
We note the decay condition on $\eps_N$ is sharp: \cite[Proposition 2.8]{alaoui2023shattering} shows that if $\eps_N\leq o(1/N)$, then the expected total variation distance between $\mu_{\beta,0}$ and $\mu_{\beta,\eps_N}$ tends to $0$ \rev{in any spherical spin glass}. Thus \eqref{eq:chaos} does not hold, i.e.\ there is no chaos.
We refer to \cite[Section 2.3]{alaoui2023shattering} for further discussion, and comparison with the overlap-based notion of disorder chaos in e.g.\ \cite{chatterjee2009disorder,chen2017parisi,chen2018energy,eldan2020simple}.

\subsection{Connections to the Kac--Rice Formula}
\label{subsec:kac-rice}

The Kac--Rice formula \cite{kac1943average,rice1944mathematical,adler2007random} forms the basis for a substantial line of work on extreme values in spherical spin glasses stemming from \cite{auffinger2013random,auffinger2013complexity}.
The formula gives general expressions for the expectations (or higher moments) of the number of critical points of $H_N$ at different energy levels (or with other characteristics), which can be studied using random matrix theory.
Notably \cite{subag2017complexity,arous2020geometry} identified the ground state energy $\GS(\xi)$ for a subclass of ``pure-like'' $\xi$, by matching the first and second moments of the number of critical points at energy $\GS_N$. 
These models exhibit $1$-RSB behavior, i.e. the minimizing function $\zeta$ in the Parisi formula is a positive constant.
Recently the author and Huang extended this argument to \emph{all} $1$-RSB models in \cite{huang2023constructive}, by truncating the second moment based on interpolation bounds.

Whenever the Kac--Rice formula suffices to identify the ground state, the radial derivative (and hence bulk Hessian spectrum) at approximate ground states can be read off from the calculation, and the Hessians will have no upward outliers in the sense of Theorem~\ref{thm:no-upward-outliers-even-xi}; \rev{see for instance \cite{xu2024hessian}.}
The reason is that the expected number of local maxima at the ground state energy level \rev{with incorrect} radial derivative (or with upward outlier eigenvalues) is exponentially small.
A similar phenomenon holds for \emph{topologically trivial} spin glasses: those whose external fields are so strong that $H_N$ has only two critical points on $\cS_N$, which are the global extrema \cite{belius2022triviality,xu2022hessian,huang2023strong}.
Topological trivialization is equivalent to the minimizing $\zeta$ being the constant function $0$.
In general, $\zeta$ can be much more complicated; see \cite{auffinger2019existence,auffinger2022spherical} for examples.
\rev{However except for the special cases of} $1$-RSB and topologically trivial models, the annealed predictions from the Kac--Rice formula are incorrect, and do not suffice to understand the behavior near the ground state 
(though see \cite{subag2017geometry,arous2024shattering} for applications to the geometry of positive temperature Gibbs measures).
In these cases, physicists have used the replica method to predict the correct quenched critical point counts; see the recent work \cite{fyodorov2018hessian,kent2023count} as well as \cite{bray1980metastable,crisanti1995thouless}.
These arguments seem difficult to make rigorous at present.
Despite these challenges, it is still natural to believe that exponentially rare behavior does not occur at approximate ground states, because the sets $A_{\delta}$ must be rather small \rev{by general principles of Gaussian process theory}.\footnote{
\rev{For instance, \cite[Theorem 1.1]{chatterjee2013level} with $\alpha=1-\frac{\delta}{\GS(\xi)}$ shows that an independent copy $H_N'$ of $H_N$ has typical ground state energy proportional to $\sqrt{\delta}$ when restricted to $A_{\delta}$.
Via Sudakov minorization, this implies $e^{o_{\delta}(N)}$ upper bounds on the typical $\eps\sqrt{N}$-covering number of $A_{\delta}$, if $\delta$ is small compared to $\eps$.
This intuitively suggests that $A_{\delta}$ is too small for other exponentially rare events to occur inside it (but does not constitute a proof since $A_{\delta}$ already depends on $H_N$).
}
}
Theorem~\ref{thm:no-upward-outliers-even-xi} confirms this intuition for the absence of upward Hessian outliers \rev{when $\xi$ is even}, but uses a completely different proof technique.

The conclusion \eqref{eq:bulk-edge-formula} is also related to thresholds obtained in \cite{auffinger2013random,auffinger2013complexity} through the Kac--Rice formula. 
Namely in pure models $\xi(t)=t^p$, for marginal stability to hold at a critical point $\bsig\in\cS_N$, one must have $H_N(\bsig)/N\approx E_{\infty}(p)=2\sqrt{\frac{p-1}{p}}$.
In mixed spherical spin glasses, a more involved Kac-Rice calculation shows that for marginally stable critical points to exist at energy $E$, one must have $E\in [E_{\infty}^-,E_{\infty}^+]$ for
\[
    E_{\infty}^{\pm}(\xi)
    \equiv
    \frac{2 \xi' \sqrt{\xi''} \pm \sqrt{4 \xi'' (\xi')^2-\left(\xi''+\xi'\right)\left(2\left(\xi''-\xi'+(\xi')^2\right)-\alpha^2 \log \frac{\xi''}{\xi'}\right)}}{\xi'+\xi''}
\]
where $\alpha=\sqrt{\xi''+\xi'-(\xi')^2}$ and $\xi,\xi',\xi''$ are evaluated at $1$.
Of course, the exact ground state is a critical point of $H_N$. 
When $\xi$ has full RSB endpoint behavior, Proposition~\ref{prop:main} implies it is marginally stable.
Thus one obtains the following inequality, where the full RSB condition is necessary because $\GS(\xi)>E_{\infty}$ for pure models. 
By contrast $\GS(\xi)\geq E_{\infty}^-(\xi)$ holds for all $\xi$; see \cite{auffinger2013complexity}, or \cite[Section 7]{huang2023strong} for extensions.

\begin{corollary}
\label{cor:E-infty}
    If $\xi$ has full RSB endpoint behavior, then $\GS(\xi)\leq E_{\infty}^+(\xi)$.
\end{corollary}

\subsection{On Marginal Stability}

The type of qualitative connection between marginal stability and full RSB proved in Corollary~\ref{cor:marginal-stability} was long anticipated in the physics literature. 
For example, the Gardner transition \cite{gardner1985spin} from $1$-RSB to full RSB in Ising spin glasses is said to occur when the Hessian spectrum at ancestor states touches $0$.
We also quote from \cite{muller2004glass}:
\[
 \textit{Here, we exploit only the well-known fact that
a full RSB glass is in a marginally stable state at all } T<T_c.
\]
See also \cite{muller2015marginal,folena2022marginal}.
It is additionally believed that ``reasonable'' optimization algorithms in high-dimension rapidly reach and then get stuck in the ``manifold'' of marginal states \cite{cugliandolo1994out,kent2024algorithm}; indeed \cite{behrens2022dis} recently conjectured that finding stable (i.e.\ non-marginal) local optima is computationally intractable in many disordered systems.
As rigorous evidence for this, \cite{huang2024optimization} analyzed the optimal message-passing algorithms to optimize $H_N$ over $\cS_N$ and showed the resulting outputs are marginally stable outside of the topologically trivial phase, even in the more general multi-species setting.
Additionally, optimal stable algorithms to optimize mean-field spin glass Hamiltonians are now understood to be closely related to full RSB \cite{subag2018following,montanari2021optimization,ams20,huang2021tight,auffinger2023optimization,huang2023algorithmic,montanari2024exceptional,jekel2024potential}.
Corollary~\ref{cor:marginal-stability} gives another rigorous link between full RSB and marginal stability, which however concerns low temperature statics rather than efficient algorithms.

\subsection{Proof Ideas}
\label{subsec:proof-ideas}

The radial derivative formula \eqref{eq:radial-deriv-formula} is proved using perturbative arguments in $\bbR^N$.
The idea is to consider the maximum value of $H_N$ on a dilated spherical domain $t\cS_N$.
For $|t-1|$ small (but independent of $N$), the change in the maximum of $H_N$ is essentially given by the maximum or minimum radial derivative value at any approximate ground state (depending on the sign of $t-1$).
On the other hand, it follows from \eqref{eq:covariance-formula} that the asymptotic ground state energy on a dilated sphere is $\GS(\xi_t)$ for some slightly perturbed $\xi_t\approx \xi$.
In this way, our proof of \eqref{eq:radial-deriv-formula} first shows an abstract formula for the radial derivative at any approximate ground state; we then substantially simplify this formula using the stationarity conditions for the minimizing $(\zeta,L)$, yielding \eqref{eq:radial-deriv-formula} and hence also \eqref{eq:bulk-edge-formula}.

The proof of \eqref{eq:top-edge-formula}, which ensures the absence of upward Hessian outliers, is more delicate and is given in Subsection~\ref{subsec:no-upward-outliers}.
Here we perturb $H_N$ by considering an augmented Hamiltonian of two replicas with constrained overlap very close to $1$.
The near-extrema for the augmented system \rev{that one should keep in mind} are obtained by starting from $\bsig\in A_{\delta}$, and perturbing $\bsig$ in opposite directions along a Hessian eigenvector.
Thus a precise estimate on the augmented ground state energy translates to a bound on the largest eigenvalue at any approximate ground state. 
We find explicit interpolation parameters in the corresponding vector Parisi formula, and perform a first-order Taylor expansion of the resulting bound for the augmented ground state energy.
This yields an eigenvalue estimate which matches the upper edge of the bulk spectrum from \eqref{eq:bulk-edge-formula}.

Finally Section~\ref{sec:alternate-proof} gives an alternate proof that typical Gibbs samples are marginally stable if $\beta\to\infty$ slowly with $N$, provided $\xi$ is generic, $1$ is not isolated in $T$, and $q\mapsto\xi''(q)^{-1/2}$ is non-affine on $(0,1]$. 
\rev{This is conceptually similar to (but strictly weaker than) Corollary~\ref{cor:marginal-stability}.}
Unlike the main proof described above, this alternate proof does not require any computations with the Parisi formula.
Instead the full RSB condition is used to ensure that low-temperature Gibbs samples have positive probability to form certain ultrametric constellations with overlap close to $1$. 
Existence of such constellations then forces the Hessian to have many near-zero eigenvalues.
This provides a more conceptual explanation for the link between marginal stability and full RSB.

\rev{We finally mention that the recent result \cite[Theorem 1.4]{dembo2024disordered} (which appeared slightly after the present work) gives an essentially complete description of the landscape around a typical Gibbs sample in spherical spin glasses, assuming genericity and an additional finite-RSB condition.
Their approach can be used to recover weaker analogs of our main results for such $\xi$ (also specialized to low temperature Gibbs samples and without a quantitative form of ``high probability'').
It is also worth mentioning that although the finite-RSB condition in \cite{dembo2024disordered} is fairly broad, it is in direct contrast with full-RSB and hence marginal stability.
}

\subsection{Technical Preliminaries}
\label{subsec:prelims}

Here we explicitly define the relevant derivative operations on the sphere, and state some useful concentration estimates.
First, the rescaled radial derivative at $\bsig\in\cS_N$ is 
\[
\partial_{\rd} H_N(\bsig) = R(\bsig, \nabla H_N(\bsig)).
\]
Next for each $\bsig \in S_N$, let $\{e_1(\bsig),\ldots,e_N(\bsig)\}$ be an orthonormal basis of $\bbR^N$ with $e_1(\bsig) = \bsig / \sqrt{N}$.
Let $\cT = \{2,\ldots,N\}$.
Let $\nabla_\cT H_N(\bsig) \in \bbR^\cT$ denote the projection of $\nabla H_N(\bsig) \in \bbR^N$ to the space spanned by $\{e_2(\bsig),\ldots,e_N(\bsig)\}$, and $\nabla^2_{\cT \times \cT} H_N(\bsig) \in \bbR^{\cT \times \cT}$ analogously. 
The spherical gradient and Hessian are defined by:
\[
    \nabla_\sph H_N(\bsig) = \nabla_{\cT} H_N(\bsig),
    \quad\quad 
    \nabla^2_{\sph} H_N(\bsig) = \nabla^2_{\cT \times \cT} H_N(\bsig) - \partial_{\rd} H_N(\bsig) I_{\cT \times \cT}.
\]
The next proposition provides smoothness estimates for $H_N$, ensuring for example that the radial derivative has typical order $O(1)$, while the spherical gradient has typical norm $O(\sqrt{N})$.
The operator norm of a tensor $\bA \in (\bbR^N)^{\otimes p}$ is
\[
    \tnorm{\bA}_\op =
    \max_{\|\bsig^1\|_2,\ldots,\|\bsig^p\|_2\leq 1}
    |\la \bA, \bsig^1 \otimes \cdots \otimes \bsig^p \ra|\,.
\]
We denote by $\sH_N\simeq\bigoplus_{p\geq 1}\bbR^{N^p}$ the set of mixed $p$-spin Hamiltonians on $\cS_N$, which for fixed $\xi$ can be identified with the coefficient \rev{tensors $\bG^{(p)}\in (\bbR^N)^{\otimes p}$ }.

\begin{proposition}
\label{prop:gradients-bounded}
    For fixed $\xi$ there exist constants $C,c>0$, and a sequence $(K_N)_{N\geq 1}$ of sets $K_N\subseteq \sH_N$, with:
    \begin{enumerate}
        \item $\bbP[H_N\in K_N]\geq 1-e^{-cN}$.
        \item If $H_N\in K_N$ and $\bx, \by\in \cS_N$, then
       \begin{align}
            \label{eq:gradient-bounded}
            \norm{\nabla^k H_N(\bx)}_{\op}
            &\le 
            CN^{1-\frac{k}{2}},\quad\forall~0\leq k\leq 3 \\
            \label{eq:gradient-lipschitz}
            \norm{\nabla^k H_N(\bx) - \nabla^k H_N(\by)}_{\op}
            &\le 
            CN^{\frac{1-k}{2}}\norm{\bx - \by},
            \quad\forall~0\leq k\leq 2.
        \end{align}
        \item If $H_N\in K_N$ then 
        \begin{equation}
        \label{eq:parisi-formula-accurate}
            \lt|\GS_N-\GS(\xi)\rt|\leq \eps.
        \end{equation}
    \end{enumerate}
\end{proposition}

\begin{proof}
    \cite[Proposition \rev{2.3}]{huang2021tight} shows that \eqref{eq:gradient-bounded} and \eqref{eq:gradient-lipschitz} hold with probability $1-e^{-cN}$. The same holds for \eqref{eq:parisi-formula-accurate} by the Borell-TIS inequality.
\end{proof}

\rev{The following description of the near-global maxima is elementary.
We explain it in detail, in part as a warm-up for the proof of Theorem~\ref{thm:no-upward-outliers-even-xi} in Subsection~\ref{subsec:no-upward-outliers}, which ends with a similar spherical Taylor expansion.
}

\begin{corollary}
    \label{cor:near-ground-state}
    Let $\delta$ be small depending on $(\xi,\eps)$.
    Then with probability $1-e^{-cN}$, the following hold for all $\bsig\in A_{\delta}$:
    \begin{align}
        \label{eq:near-ground-state}
        |H_N(\bsig)/N-\GS(\xi)|&\leq \eps,
        \\
        \label{eq:near-critical-point}
        \|\nabla_{\sph} H_N(\bsig)\|&\leq \eps\sqrt{N},
        \\
        \label{eq:near-concave}
        \blambda_1\big(\nabla_{\sph}^2 H_N(\bsig)\big)&\leq \eps.
    \end{align}
\end{corollary}

\begin{proof}
    The first bound \eqref{eq:near-ground-state} follows by \eqref{eq:parisi-formula-accurate}.
    \rev{
    To deduce \eqref{eq:near-critical-point} and \eqref{eq:near-concave}, we will show that if one of them does not hold, then for $H_N\in K_N$, then Taylor expanding $H_N$ near $\bsig$ would falsify \eqref{eq:parisi-formula-accurate}.
    In both cases, it will be convenient to use the following spherical Taylor expansion, which is easily verified using Proposition~\ref{prop:gradients-bounded} and justifies the definition of the Riemannian Hessian $\nabla_{\sph}^2$.
    For $\bv\in \cS_N\cap \bsig^{\perp}$, we have $\cos(\eps)\bsig + \sin(\eps)\bv\in \cS_N$ and:
    \begin{equation}
    \label{eq:taylor-expansion-HN}
    \begin{aligned}
    H_N(\cos(\eps)\bsig + \sin(\eps)\bv)
    &=
    H_N(\bsig)
    +
    \eps \la \nabla_{\sph} H_N(\bsig),\bv\ra
    +
    \frac{\eps^2}{2}
    \la \nabla_{\sph}^2 H_N(\bsig),\bv^{\otimes 2}\ra
    +
    O(\eps^3 N)
    \\
    &=
    H_N(\bsig)
    +
    \eps \la \nabla_{\sph} H_N(\bsig),\bv\ra
    +
    O(\eps^2 N).
    \end{aligned}
    \end{equation}
    Here the hidden implicit constants are uniform in $\bsig$ and $\bv$ and $\eps$, for all $H_N\in K_N$.
    }

    \rev{
    For \eqref{eq:near-critical-point}, let us assume $H_N\in K_N$ and (without loss of generality) that $\eps$ is small depending on $\xi$.
    Suppose that $\delta<\eps^3/4$ is such that \eqref{eq:parisi-formula-accurate} holds with constant $\eps^3/4$ (noting that $\eps>0$ was arbitrary in the original statement and hence can be adjusted separately in the three estimates to be proved).
    We will show that \eqref{eq:near-critical-point} then follows with constant $\eps$.
    Indeed if it does not, we let $\bv=\sqrt{N}\nabla_{\sph}H_N(\bsig)/\|\nabla_{\sph}H_N(\bsig)\|$ and 
    $\eps'=\eps \|\nabla_{\sph}H_N(\bsig)\|/\sqrt{N}$.
    Then 
    \begin{align*}
    H_N(\cos(\eps')\bsig+\sin(\eps')\bv)
    &=
    H_N(\bsig)
    +
    \eps\|\nabla_{\sph}H_N(\bsig)\|^2
    +
    O(\eps'^2 N)
    \\
    &\geq 
    H_N(\bsig)
    +
    \frac{\eps\|\nabla_{\sph}H_N(\bsig)\|^2}{2}
    \\
    &\geq 
    H_N(\bsig)
    +
    \frac{\eps^3 N}{2}
    \\
    &> 
    N\Big(\GS(\xi)+\frac{\eps^3}{4}\Big)
    .    
    \end{align*}
    Here the second inequality follows from the assumption that $\|\nabla_{\sph}H_N(\bsig)\|\geq \eps\sqrt{N}$ and we used $\bsig\in A_{\delta}$ for $\delta<\eps^3/4$ to obtain the last inequality.
    Comparing the initial and final expressions above contradicts our assumption that \eqref{eq:parisi-formula-accurate} holds with constant $\eps^3/4$.

    For the last assertion, we again suppose $H_N\in K_N$ and $\eps$ is small depending on $\xi$.
    Here we assume \eqref{eq:near-concave} does not hold with constant $\eps$ but that both \eqref{eq:near-critical-point} and \eqref{eq:parisi-formula-accurate} hold with constant $\eps^5/8$.
    Then let $\bw/\sqrt{N}$ be the top unit eigenvector of $\nabla_{\sph}^2 H_N(\bsig)$, viewed as a quadratic form on $\bsig^{\perp}$.
    With $\eps''=\eps \blambda_1$ for $\blambda_1=\blambda_1(\nabla_{\sph}^2 H_N(\bsig))>\eps$ we find 
    \begin{align*}
    H_N(\cos(\eps'')\bsig+\sin(\eps'')\bw)
    &=
    H_N(\bsig)
    +
    O(\eps'' \|\nabla_{\sph}H_N(\bsig)\|\sqrt{N})
    +
    \frac{(\eps'')^2 \blambda_1 N}{2}
    +
    O((\eps'')^3 N)
    \\
    &\geq 
    H_N(\bsig)
    +
    (\eps'')^2\blambda_1 N/4
    \\
    &\geq 
    H_N(\bsig)
    +
    \eps^5 N/4
    \\
    &> 
    N\Big(\GS(\xi)+\frac{\eps^5}{8}\Big).
    \end{align*}
    Again we used $\bsig\in A_{\delta}$ (now with $\delta<\eps^5/8$) to obtain the last inequality, and the final inequality contradicts the assumption that \eqref{eq:parisi-formula-accurate} holds (with constant $\eps^5/8$).
    } 
\end{proof}

The next standard lemma states that $\partial_{\rd}H_N(\bsig)$ determines the bulk spectrum of $\nabla_{\sph}^2 H_N(\bsig)$ uniformly over $\bsig\in\cS_N$.
In particular, \eqref{eq:bulk-edge-formula} and \eqref{eq:radial-deriv-formula} are equivalent (the case $N(1-\eta)\leq j\leq N-K$ follows by negating $H_N$, which preserves its law).
We note that in Section~\ref{sec:alternate-proof}, it is important to take $j$ constant in Lemma~\ref{lem:spectrum-approx}.

\begin{lemma}[{\cite[\rev{Lemma 1.3}]{subag2018following}}]
\label{lem:spectrum-approx}
    For any $\eps>0$ there are $K=K(\xi,\eps)$ and $\eta=\eta(\xi,\eps)>0$ such that
    \[
    \big|
    \blambda_j\big(\nabla_{\sph}^2 H_N(\bsig)\big)
    -
    2\xi''(1)^{\rev{1/2}}
    \rev{+}
    \partial_{\rd}H_N(\bsig)
    \big|
    \leq
    \eps
    \]
    holds for all $K\leq j\leq \eta N$ and $\bsig\in\cS_N$ simultaneously, with probability $1-e^{-cN}$.
\end{lemma}

Finally we record the standard fact that low temperature Gibbs samples are approximate ground states.

\begin{lemma}
\label{lem:small-values-low-prob}
    Assume $H_N\in K_N$. Then for $\beta\geq \beta_*(\xi)$ large, $A_{\beta^{-1/2}}$ (recall \eqref{eq:approx-ground-states}) satisfies $\mu_{\beta}(A_{\beta^{-1/2}})\geq 1-e^{-cN}$.
\end{lemma}

\begin{proof}
    Let $\bsig_*\in\cS_N$ be the maximizer of $H_N$.
    The radius $\beta^{-2}\sqrt{N}$ neighborhood of $\bsig_*$ has $\mu_0$ (uniform) measure $\beta^{-O(N)}$ and energy within $O(\beta^{-2}N)$ of the maximum. The contribution to $Z_{N,\beta}$ from this neighborhood is
    \[
    \exp\Big(\beta H_N(\bsig_*) - N\cdot O(\log\beta)\Big).
    \]
    When $\beta$ is large, this is exponentially larger than the contribution to $Z_{N,\beta}$ from the complement of $A_{\beta^{-1/2}}$.
\end{proof}

\subsection{Characterization of the Minimizer in the Parisi Formula}
\label{subsec:characterization-of-minimizer}

Since the Parisi functional $\cQ$ is strictly convex, the unique minimizer $(\zeta,L)$ is characterized by first-order stationarity conditions. 
These will be useful later, and we review them now.
Given $(\zeta,L)\in\cK$ (recall \eqref{eq:def-cK}), define
\begin{align}
    \label{eq:def-G}
    G(q) &= \xi'(q) - \int_0^q \fr{\de s}{\hat \zeta(s)^2}, &
    g(s) &= \int_s^1 G(q)~\de q.
\end{align}
Let $\nu$ be the finite Borel measure on $[0,1]$ defined by 
\begin{equation}
    \label{eq:def-nu-infty}
    \nu([0,q]) = \zeta(q) \qquad \forall q\in [0,1]
\end{equation}
and define the set
\begin{equation}
\label{eq:T-def}
    T = \{q \in [0,1] : g(q) = 0\}.
\end{equation}
Note that $1\in T$ holds trivially.

The characterization below is primarily from \cite[Theorem 2]{chen2017parisi}. 
In the finite temperature case \cite[Corollary 1.6]{jagannath2018bounds} proves that $T$ consists of finitely many intervals, and the zero temperature case is identical. (One combines \cite[Theorem 1.13]{jagannath2017low} and the observation that by analyticity, $\big(\frac{1}{\sqrt{\xi''}}\big)''$ changes sign finitely many times on $[0,1]$.)

\begin{proposition}
    \label{prop:cs-extremality-0temp}
    The unique $(L,\zeta)$ attaining the infimum \eqref{eq:cs-functional-0temp-inf} is characterized as follows, where $\supp(\zeta)\subseteq [0,1)$ denotes the set of points of increase of $\zeta$ (i.e.\ the support of $\nu$ in \eqref{eq:def-nu-infty}):
    \[
        G(1)=0;\quad\quad
        \min_{q\in [0,1]} g(q) = 0;\quad\quad
        T\rev{\supseteq} \supp(\zeta).
    \]
    Furthermore, $T$ is a disjoint union of finitely many closed intervals (possibly including singletons).
\end{proposition}

\begin{corollary}
\label{cor:full-RSB-justification}
    If $1$ is not isolated in $T$, then $\xi$ exhibits full RSB endpoint behavior.
    If $q\mapsto\xi''(q)^{-1/2}$ is non-affine on $(0,1]$, then $\zeta$ is strictly increasing near $1$.
\end{corollary}

\begin{proof}
    Proposition~\ref{prop:cs-extremality-0temp} implies that $T$ contains a non-trivial interval $(a,1)$.
    Set $f(q)=\xi''(q)^{-1/2}$. On this interval, we must have $g''(q)=0$, which gives
    \[
    \hat\zeta(q)=f(q),
    \qquad
    \zeta(q)=-f'(q),
    \qquad
    \zeta'(q)=-f''(q).
    \]
    Evaluating the first identity at $q=1$ yields full RSB endpoint behavior.
    For the latter assertion, analyticity on the connected interval $(0,1]$ and the non-affineness hypothesis imply that $f''$ is not identically zero on $(a,1)$.
    Moreover, $f''$ is analytic on an open complex neighborhood of $[a,1]$.
    Thus it has only finitely many zeros on a sufficiently small interval ending at $1$.
    Since $\zeta$ is non-decreasing, $-f''=\zeta'$ is non-negative there, and hence positive away from those finitely many zeros.
    It follows that $\zeta$ is strictly increasing on $(1-\eps_*(\xi),1)$.
\end{proof}

\begin{remark}
\label{rem:affine-exception}
The non-affineness hypothesis in the second assertion is necessary.
For example, for $\kappa>0$ and $0<\rho<1$, let
\[
\xi(q)=\frac{1}{\kappa^2}\bigl(-\log(1-\rho q)-\rho q\bigr).
\]
This is a valid mixture function because
\[
\xi(q)=\frac{1}{\kappa^2}\sum_{p=2}^{\infty}\frac{\rho^p}{p}q^p.
\]
Thus its coefficients $\gamma_p^2=\rho^p/(p\kappa^2)$ are non-negative, and its radius of convergence is $1/\rho>1$.
Then
\[
\xi''(q)^{-1/2}=\frac{\kappa}{\rho}-\kappa q,
\]
and the Parisi minimizer is $\zeta(q)\equiv\kappa$ with $L=\kappa/\rho$.
Thus $T=[0,1]$ and the full RSB endpoint equality holds, even though $\zeta$ is constant near $1$.
\end{remark}

\revv{
\begin{remark}
\label{rem:not-equivalence}
The converse of the first assertion in Corollary~\ref{cor:full-RSB-justification} need not hold.
Indeed, suppose the non-trivial interval $[t,s]\subseteq (0,1)$ is a connected component of the set $T$ (see \cite{auffinger2022spherical} for explicit constructions of such behavior).
As in the proof of Corollary~\ref{cor:full-RSB-justification}, one then finds $\hat\zeta(t)=\xi'(t)^{-1/2}$.
It is not hard to show that $\xi_{\sqrt{t}}(q)=\xi(tq)$ has full-RSB endpoint behavior at $1$, but that $1$ is an isolated point of $T(\xi_{\sqrt{t}})$.
This is because the associated Parisi order parameter is a reparametrization of $\zeta$ (see e.g. \cite[Proposition 11]{subag2018free} or \cite[Proposition 4.3]{huang2023constructive}).
We expect such counterexamples to essentially comprise the boundary of the full RSB phase.
\end{remark}
}

\revv{
\begin{remark}
\label{rem:weaker-condition}
As previously mentioned, a large portion of Corollary~\ref{cor:deep-level-sets} follows from the weaker assumption that \emph{strict} full RSB endpoint behavior is absent (i.e.\ that $1$ is isolated in $T$), without requiring the approach in Theorem~\ref{thm:no-upward-outliers-even-xi}.
We assume $(a,1)$ is disjoint from $T$, and again that $\xi$ is even.
It follows from two-replica interpolation bounds at zero temperature (see e.g. \cite[Theorem 6 and 11]{auffinger2018energy}, or \cite[Theorem 3]{chen2019suboptimality} as summarized in \cite[Lemma 6.4]{sellke2021optimizing}) that for small $\delta$ depending on $\eps$, with probability $1-e^{-cN}$ there does not exist any pair $(\bsig,\bsig')\in A_{\delta}\times A_{\delta}$ such that $R(\bsig,\bsig')\in (a+\eps,1-\eps)$.
Taking $\eps$ small depending on $a$, we may define a graph on vertex set $A_{\delta}$ such that $\bsig,\bsig'$ are adjacent if $R(\bsig,\bsig')\geq 1-\eps$.
Then the property above implies all connected components of this graph have Euclidean diameter at most $2\sqrt{\eps N}$, and are pairwise $10\sqrt{\eps N}$-separated.

To our knowledge these components need not be connected with respect to the usual topology on $\cS_N$, and indeed could have many critical points for $H_N$.
Still, the existence of such a clustering suffices for our proofs of Corollaries~\ref{cor:langevin} and \ref{cor:disorder-chaos} to go through unchanged.
We thank the anonymous referee for helpful observations that led to this remark.
\end{remark}
}

\section{Proofs of Main Results}
\label{sec:proofs}

Here we prove Proposition~\ref{prop:main}, Theorem~\ref{thm:no-upward-outliers-even-xi}, and Corollary~\ref{cor:disorder-chaos}.
Subsection~\ref{subsec:radial-derivative} computes the radial derivative at approximate ground states, thus proving \eqref{eq:radial-deriv-formula}.
As discussed previously, this is done by considering the maximum value of $H_N$ on slight dilations of $\cS_N$, and then simplifying the resulting abstract formula using the stationarity conditions from Subsection~\ref{subsec:characterization-of-minimizer}.
The next subsection then proves \eqref{eq:top-edge-formula} using a two-replica interpolation bound.
Finally Subsection~\ref{subsec:chaos} explains how to deduce Corollary~\ref{cor:disorder-chaos} from Corollary~\ref{cor:deep-level-sets}.

\subsection{Radial Derivative at Approximate Ground States}
\label{subsec:radial-derivative}

We fix $\xi$ and $\eta<0.01$ with $\xi(1+\rev{2\eta})<\infty$, and let $\xi_t(q)=\xi(t^2 q)$ for $t\in (1-\eta,1+\eta)$.
To explain this definition, note that by \eqref{eq:covariance-formula}, the function on $\cS_N$ defined by $\bsig\to H_N(t\bsig)$ is precisely a spherical spin glass with mixture $\xi_t$.

We would like to show $t\mapsto \GS(\xi_t)$ is in $C^1([1-\eta,1+\eta])$ using the envelope theorem.
A general differentiation formula with respect to each coefficient $\gamma_p$ is stated in \cite[Remark 2]{chen2017parisi}, analogously to the positive temperature case; this of course suggests a formula for $\frac{\de}{\de t}\GS(\xi_t)$ by linearity. For completeness we give a careful proof, primarily checking that the minimizers remain within a suitable compact subset of $\cK$.\footnote{To illustrate the need for care, note that $\cK$ is not compact, and $\zeta\mapsto \int_0^1 \zeta(q)\de q$ is lower semi-continuous but not continuous in the vague topology. 
Lower semi-continuity does not imply $\inf\limits_{|t-1|\leq\eta}\hat\zeta_t(1)>0$, even given that $t\mapsto (\zeta_t,L_t)\in\cK$ is continuous and $\hat\zeta_t(1)>0$ for fixed $t$. 
}
Thus, let 
\[
(\zeta_t,L_t)=\argmin_{(\zeta,L)\in\cK} \cQ(\zeta,L;\xi_t)
\]
be the corresponding \rev{(unique)} minimizers in the zero-temperature Parisi formula.

\begin{lemma}
\label{lem:compactification}
    There exists $C=C(\xi,\eta)>0$ such that for all $t\in (1-\eta,1+\eta)$,
    \begin{equation}
    \label{eq:zero-temp-parisi-nondegenerate}
    L_t\leq C;\quad \quad 
    \hat\zeta_t(1)\geq 1/C.
    \end{equation}
\end{lemma}

\begin{proof}
    We treat $C=C(\xi,\eta)$ as a constant that may vary from line to line.
    In the ``topologically trivial'' case that $\zeta_t$ is identically $0$, we find $\revv{L_t=\hat\zeta_t(1)=1/\sqrt{\xi_t'(1)}}$, so the conclusion is obvious.
    We assume this is not the case below.
    Noting that $\hat\zeta_t(\cdot)$ is a decreasing function 
    \rev{for each $t$}, it follows that $\hat\zeta_t(1/2)\leq C$ is bounded independently of $t$ since
    \begin{align*}
    \rev{2\,\GS(\xi_{1+\eta})}
    \geq 2\,\GS(\xi_t)
    &\geq
    2\cQ(\zeta_t,L_t;\xi_t)
    \\
    &\geq 
    \int_0^1 \xi_t''(q)\rev{\hat\zeta_t(q)}\de q
    \\
    &\geq 
    \int_{0}^{1/2} \xi_t''(q)\rev{\hat\zeta_t(q)}\de q
    \\
    &\geq 
    \frac{\hat\zeta_t(1/2)}{2}
    \int_0^{1/2} \xi_t''(q)~\de q
    \\
    &\geq 
    \frac{\hat\zeta_t(1/2)}{2}
    \revv{\big(\xi_t'(1/2)-\xi_t'(0)\big)}
    \\
    &\geq 
    \Omega(\hat\zeta_t(1/2)).
    \end{align*}
    \rev{(Since $\xi(1+2\eta)<\infty$, the coefficients of $\xi_{1+\eta}$ decay exponentially and so $\GS(\xi_{1+\eta})<\infty$. We also used the easy fact that $GS(\xi_t)$ is increasing in $t$ in the first step.)}
    \rev{Next we note that $\hat\zeta_t(\cdot)$ is concave and non-negative on $[0,1]$ since $\zeta_t(\cdot)$ is increasing. 
    Therefore 
    \[
    \hat\zeta_t(0)\leq 2\hat\zeta_t(1/2)-\hat\zeta_t(1)
    \leq 
    2\hat\zeta_t(1/2)
    \leq 2C
    \]
    is bounded depending only on $\xi$.
    Since $L_t=\hat\zeta_t(0)$ we conclude the first estimate.
    }

    We now turn to the second estimate.
    \rev{
    Let 
    \revv{
    \[
    A=
    \frac{1}{2\sqrt{\GS(\xi_{1+\eta})^2+\xi_{1+\eta}''(1)}}
    \leq 
    \frac{1}{2\sqrt{\GS(\xi_{t})^2+\xi_{t}''(1)}}
    \leq 
    \frac{1}{2\GS(\xi_t)}
    \]
    }
    and assume that $\hat\zeta_t(1)<A$ (else we are already done).
    We first claim that $\hat\zeta_t(0)\geq A$.
    Indeed if not, then \revv{$\hat\zeta_t(q)<\frac{1}{2\GS(\xi_t)}$ for all $q\in [0,1]$ 
    and so 
    \[
    \cQ(\zeta_t,L_t;\xi_t)
    \geq 
    \frac{1}{2}\int_0^1 \frac{\de q}{\hat\zeta_t(q)}
    >\GS(\xi_t)
    \]
    }
    which is impossible.
    Given that $\hat\zeta_t(0)\geq A>\hat\zeta_t(1)$, it follows by continuity that there exists $q_*\in [0,1)$ such that $\hat\zeta_t(q_*)=A$.
    We will show that $\zeta_t$ is constant on $[q_*,1]$.
    Indeed consider the modified order parameter 
    \[
    \rev{\zeta^*_t}(q)=
    \begin{cases}
    \zeta(q),\quad q\leq q_*;
    \\
    \zeta(q_*),\quad q\geq q_*.
    \end{cases}
    \]
    We aim to show that $\zeta^*_t(q)=\zeta_t(q)$ for all $q\in [0,1]$, \revv{via the inequality $\cQ(\zeta_t^*,L_t;\xi_t)\leq \cQ(\zeta_t,L_t;\xi_t)$}.
    Setting
    \[
    \hat\zeta^*_t(q)\equiv 
    L_t-\int_0^q \zeta^*_t(u)\de u,
    \]
    by construction $\hat\zeta_t^*$ is decreasing with $\hat\zeta_t^*(q_*)=A$,
    \revv{and so $\hat\zeta_t^*(q)\leq A$ for all $q\geq q_*$}.
    Further $\zeta^*_t(q)\leq \zeta_t(q)$ for each $q\in [0,1]$, which implies by integration that $\hat\zeta^*_t(q)\geq \hat\zeta_t(q)$ for all $q$.
    \revv{Combining, we have for $q\geq q_*$ that
    \[
    \hat\zeta_t^*(q)\in [\hat\zeta_t(q),A]
    \subseteq 
    [\hat\zeta_t(q),\xi_t''(q)^{-1/2}]
    \]
    which by monotonicity of $x\mapsto \xi_t''(q)x+\frac{1}{x}$ on $x\in (0,\xi_t''(q)^{-1/2})$ implies 
    \begin{equation}
    \label{eq:AM-GM}
    \xi_t''(q)\hat\zeta_t(q)+\frac{1}{\hat\zeta_t(q)}
    \geq 
    \xi_t''(q)\hat\zeta^*_t(q)+\frac{1}{\hat\zeta^*_t(q)}.
    \end{equation}
    }
    For $q\leq q_*$, we also have \eqref{eq:AM-GM}, in fact with equality.
    Inspecting the definition of $\cQ$, it then follows that 
    \[
    \cQ(\zeta_t,L_t;\xi_t)\ge \cQ(\zeta^*_t,L_t;\xi_t).
    \]
    By minimality of $(\zeta_t,L_t)$ we must have equality, which means equality holds in \eqref{eq:AM-GM} for every $q$.
    Hence $\zeta_t=\zeta^*_t$ and so $\zeta_t(q)$ must be constant on $q\in [q_*,1)$.
    Below we let $z=\zeta_t(q_*)$ be this constant value.
    }
    
    Next, let $q^*\in [0,q_*]$ be the largest point in the support of $\zeta_t$
    \rev{(we already ruled out the case that $\zeta_t$ is identically zero).}
    Recalling the notation from \eqref{eq:def-G}, we have $G(1)=G(q^*)$ which means 
    \[
    \xi_t'(1)-\xi_t'(q^*)
    =
    \int_{q^*}^1 \frac{\de q}{\hat\zeta_t(q)^2}
    =
    \frac{1}{z\hat\zeta_t(1)}
    -
    \frac{1}{z\hat\zeta_t(q^*)}
    .
    \]
    Note that $\hat\zeta_t(q^*)\geq \hat\zeta_t(q_*)=A$.
    If $\hat\zeta_t(1)\geq A/2$ we are done; if not, we conclude that 
    \[
    \xi_t'(1)-\xi_t'(q^*)
    \geq 
    \frac{1}{2z\hat\zeta_t(1)}
    \implies 
    \hat\zeta_t(1)\geq \frac{1}{(1-q^*)z C'(\xi)},
    \]
    \rev{where 
    \[
    C'(\xi)
    =
    2\xi_{1+\eta}''(1)
    \geq 
    2\xi_{t}''(1)
    \geq 
    2\cdot \frac{\xi_t'(1)-\xi_t'(q^*)}{1-q_*}.
    \]
    }
    However by our assumptions and the first part above, 
    \[
    (1-q^*)z
    =
    \rev{
    \int_{q^*}^1 \zeta_t(q)dq
    }
    =
    \rev{\hat\zeta_t(q^*)-\hat\zeta_t(1)}
    \leq L_t\leq C
    .
    \]
    Thus $\hat\zeta_t(1)$ is bounded away from $0$ depending only on $\xi$ in all cases.
\end{proof}

Given $C>0$, let $\cK_C\subseteq \cK$ consist of those $(\zeta,L)$ satisfying \eqref{eq:zero-temp-parisi-nondegenerate}.
We endow $\cK_C$ with the weak*-topology from the associated $\nu$ in \eqref{eq:def-nu-infty}, making it a compact metric space for each $C$.

\begin{proposition}
\label{prop:ground-state-derivative}
    For $s\in (1-\eta,1+\eta)$, the function $s\mapsto \GS(\xi_s)$ is continuously differentiable with derivative
    \begin{equation}
    \label{eq:abstract-formula}
    \rev{\frac{1}{s}}
    \frac{\de}{\de s}
    \GS(\xi_s)
    =
    \rev{\frac{1}{s}}
    \lt(\frac{\de}{\de t}
    \cQ(\zeta_s,L_s;\xi_t)\big|_{t=s}
    \rt)
    =
    \xi_s'(0)L_s
    +
    \int_0^1 
    \big(2\xi_s''(q)+q\xi_s'''(q)\big)\hat\zeta_s(q)
    ~\de q.
    \end{equation}
\end{proposition}

\begin{proof}
    Lemma~\ref{lem:compactification} shows $\GS(\xi_s)$ is defined as an infimum over the compact set $\cK_C\subseteq \cK$; Proposition~\ref{prop:cs-extremality-0temp} ensures the minimizer is unique.
    This easily implies $s\mapsto (\zeta_s,\rev{L_s})$ is continuous.
    \rev{Using the just-established compactness}, the first equality now follows from the envelope theorem \cite[Theorem 2 and Corollary 4]{milgrom2002envelope}, assuming said derivative of $\cQ$ exists.
    Indeed dominated convergence \rev{and the chain rule yields} the explicit formula above. 
    (The fact that $C$ is fixed also obviates any analytic issues near $q\approx 1$.)
\end{proof}

Next we show this formula coincides with $r(\xi)$ as defined in \eqref{eq:r-def}, using the stationarity conditions for $(\zeta,L)$.

\begin{lemma}
\label{lem:integration-by-parts-for-r-xi}
    $r(\xi_s)$ agrees with the formula \eqref{eq:abstract-formula}.
\end{lemma}

\begin{proof}
    We use the stationarity conditions reviewed in Subsection~\ref{subsec:characterization-of-minimizer}, and set $s=1$ for convenience.
    Since $\hat\zeta'(q)=-\zeta(q)$, integrating by parts gives 
    \[
    \int_0^1 
    \big(\xi''(q)+q\xi'''(q)\big)\hat\zeta(q)
    \de q
    =
    \underbrace
    {[q\xi''(q)\hat\zeta(q)]\big|_0^1}
    _
    {\xi''(1)\hat\zeta(1)}
    +
    \int_0^1 q\xi''(q)\zeta(q)\de q.
    \]
    Recalling the definition \eqref{eq:r-def} of $r(\xi)$, it remains to show that 
    \begin{equation}
    \label{eq:remains-to-show}
    \xi'(0)L
    +
    \int_0^1 
    \xi''(q)\hat\zeta(q)+q\xi''(q)\zeta(q)
    \de q
    \stackrel{?}{=} \hat\zeta(1)^{-1}.
    \end{equation}
    Note that by the last assertion in Proposition~\ref{prop:cs-extremality-0temp}, there exists a unique finite sequence 
    $0\leq q_0<q_1<\dots <q_D=1$ such that for each $0\leq d\leq D$, either one of the closed intervals $[q_{d-1},q_d]$ or $[q_d,q_{d+1}]$ is a connected component of $T$, or $q_d$ is an isolated point.    
    We verify \eqref{eq:remains-to-show} by evaluating the integral separately over each subinterval $[q_d,q_{d+1}]$.

    We first handle the contribution from non-trivial intervals $[q_d,q_{d+1}]\subseteq T$.
    In this case, \eqref{eq:def-G} gives $g''(q)=0$ for all $q\in [q_d,q_{d+1}]$, implying $\hat\zeta(q)=\xi''(q)^{-1/2}$. 
    Then $\zeta(q)=\frac{\xi'''(q)}{2\xi''(q)^{3/2}}$, so the contribution to the integral in \eqref{eq:remains-to-show} is:
    \[
    \int_{q_d}^{q_{d+1}}
    \xi''(q)^{1/2}+\frac{q\xi'''(q)}{2\xi''(q)^{1/2}}\de q
    =
    [q\xi''(q)^{1/2}]\big|_{q_d}^{q_{d+1}}
    =
    \frac{q_{d+1}}{\hat\zeta(q_{d+1})}
    -
    \frac{q_{d}}{\hat\zeta(q_{d})}.
    \]
    
    Next suppose $[q_d,q_{d+1}]\cap T=\{q_d,q_{d+1}\}$.
    Then $\zeta$ is constant on $(q_d,q_{d+1})$, and so on this interval $\hat\zeta(q)=x-yq$ for some real $(x,y)=(x_d(\xi),y_d(\xi))$ with $y\geq 0$, where $\zeta(q)=y$.
    Then the contribution is 
    \[
    \int_{q_d}^{q_{d+1}} 
    \xi''(q)(x-yq)+q\xi''(q)y\de q
    =
    x\big(\xi'(q_{d+1})-\xi'(q_d)\big).
    \]
    Since $g$ is minimized at both $q_d$ and $q_{d+1}$, we have $G(q_d)=G(q_{d+1})=0$.
    Thus we similarly find:
    \begin{align*}
    x\big(\xi'(q_{d+1})-\xi'(q_d)\big)
    &=
    x\int_{q_d}^{q_{d+1}} \frac{\de q}{\hat\zeta(q)^2}
    =
    \int_{q_d}^{q_{d+1}} \frac{x}{(x-yq)^2} \de q
    \\
    &=
    \lt[\frac{q}{x-yq}\rt]\Big|_{q_d}^{q_{d+1}}
    =
    \frac{q_{d+1}}{\hat\zeta(q_{d+1})}
    -
    \frac{q_{d}}{\hat\zeta(q_{d})}
    .
    \end{align*}
    Telescoping,
    \[
    \int_{q_0}^{1}
    \xi''(q)\hat\zeta(q)+q\xi''(q)\zeta(q)
    \de q
    =
    \frac{1}{\hat\zeta(1)}-\frac{q_0}{\hat\zeta(q_0)}.
    \]
    If $q_0=0$, this completes the proof (this condition is equivalent to $\xi'(0)=0$).
    If not, since $\zeta(q)=0$ for $q<q_0$, the function $\hat\zeta$ is constant on $[0,q_0]$ and so 
    \[
    \int_0^{q_0} \xi''(q)\hat\zeta(q)~\de q
    =
    L\big(\xi'(q_0)-\xi'(0)\big).
    \]
    Furthermore $G(q_0)=0$, since either $q_0\in (0,1)$ must be an interior minimizer of $g$ (which implies $G(q_0)=-g'(q_0)=0$), or $q_0=1$ so again $G(q_0)=0$ by Proposition~\ref{prop:cs-extremality-0temp}.
    Thus
    \[
    \xi'(q_0)
    =
    \int_0^{q_0}\frac{\de q}{\hat\zeta(q)^2}
    =
    \frac{q_0}{\hat\zeta(0)^2}
    =
    \frac{q_0}{L^2}
    .
    \]
    Thus the first term in \eqref{eq:remains-to-show} is $\xi'(0)L=\frac{q_0}{L}=\frac{q_0}{\hat\zeta(q_0)}$ when $q_0>0$.
    Matching terms finishes the proof.
\end{proof}

Next we apply Proposition~\ref{prop:ground-state-derivative} and standard concentration estimates to determine the radial derivative at any near-ground state, thus proving \eqref{eq:radial-deriv-formula}.
The point is that, as explained at the start of this subsection, if $\xi$ is the mixture function for $H_N(\bsig)$, then $\xi_t$ is the mixture function for $H_N(t\bsig)$.
On the other hand, elementary calculus in $\bbR^N$ shows the radial derivative at ground states determines the change in the maximum value of $H_N$ from a small dilation of $\cS_N$.
(See \cite{auffinger2018concentration} for general related results.)

\begin{proof}[Proof of Proposition~\ref{prop:main}]
    Due to Lemma~\ref{lem:spectrum-approx}, it suffices to prove \eqref{eq:radial-deriv-formula}.
    In light of \eqref{eq:near-critical-point} in Corollary~\ref{cor:near-ground-state}, the nontrivial portion of \eqref{eq:radial-deriv-formula} is to compute the radial derivative.
    We restrict to the event $K_N$ of Proposition~\ref{prop:gradients-bounded}.
    Let $H_N$ be a random Hamiltonian with mixture $\xi$ and let
    \[
    \GS(H_N; t)
    =
    \max_{\bsig\in\cS_N} H_N(t\bsig)/N.
    \]
    Setting $\eta=\sqrt{\delta}$, we consider the event that
    \begin{equation}
    \label{eq:event-GS-typical}
    |\GS(H_N; t)-\GS(\xi_t)|\leq \delta,
    \quad\quad
    \forall\, t\in \{1-\eta,1,1+\eta\}
    .
    \end{equation}
    It follows from the discussion above, e.g.\ \eqref{eq:covariance-formula}, that this event has probability at least $1-e^{-cN}$.
    We will show it implies the desired conclusion.
    
    Assume for sake of contradiction that $\partial_{\rd}H_N(\bsig)\geq r(\xi)+\eps$ for some $\delta$-approximate ground state $\bsig$; the opposite case is similar using $1-\eta$ instead of $1+\eta$.
    Taylor expanding and bounding the error terms with Proposition~\ref{prop:gradients-bounded}, 
    \[
    \frac{H_N((1+\eta)\bsig)-H_N(\bsig)}{N}
    \geq 
    \eta \partial_{\rd}H_N(\bsig)
    -
    C(\xi)\eta^2.
    \]
    Since $\bsig$ is an $\eta^2$-approximate ground state for $H_N$, we find from \eqref{eq:event-GS-typical} and the assumption on $\partial_{\rd}H_N(\bsig)$ that
    \begin{align*}
    \GS(\xi_{1+\eta})
    &\geq 
    H_N((1+\eta)\bsig)/N - \delta
    \\
    &\geq 
    H_N(\bsig)/N + \eta\partial_{\rd}H_N(\bsig)-(C(\xi)+1)\eta^2
    \\
    &\geq 
    \GS(\xi_1)
    +
    \eta \partial_{\rd}H_N(\bsig)
    -
    (C(\xi)+2)\eta^2
    \\
    &\geq 
    \GS(\xi_1)
    +
    \eta r(\xi)
    +
    \eta\eps
    -
    (C(\xi)+2)\eta^2
    .
    \end{align*}
    On the other hand Proposition~\ref{prop:ground-state-derivative} ensures the continuous differentiability of $\GS(\xi_t)$, so for $\eta\leq \eta_*(\xi,\eps)$ small,
    \[
    \GS(\xi_{1+\eta})
    \leq 
    \GS(\xi_1)
    +
    \eta r(\xi)
    +
    \frac{\eta\eps}{2}.
    \]
    The previous two displays are contradictory for $\eta$ small depending on $(\xi,\eps)$, which completes the proof.
\end{proof}

\subsection{No Upward Outliers in Even Models}
\label{subsec:no-upward-outliers}

When $\xi(q)=\sum_{p\in 2\bbN} \gamma_p^2 q^p$ is even, we show that the Hessian at near ground states has no eigenvalues above the upper edge of the bulk spectrum. 
We prove \eqref{eq:top-edge-formula}, the absence of eigenvalues above the bulk, by estimating the ground state energy of a two-replica spin glass:
\begin{equation}
\label{eq:GS-2-eps}
\GS_{2,\eps}(H_N)
=
\sup_{\substack{\bsig,\bsig'\in\cS_N\\ R(\bsig,\bsig')=1-\eps}}
\frac{H_N(\bsig)+H_N(\bsig')}{N}.
\end{equation}
When $\xi$ is even, the limit $\GS_{2,\eps}(\xi)=\plim_{N\to\infty} \GS_{2,\eps}(H_N)$ exists and is given by a two-dimensional constrained generalization of the Parisi formula.
Although we refrain from giving the notationally heavy general definitions, vector models of this type have played an important role in spin glass theory, being used in Talagrand's original proof of the Parisi formula \cite{talagrand2006parisi,talagrand2006free} and subsequent works \cite{panchenko2007overlap,chen2015disorder,chen2017variational,chen2018energy,auffinger2018energy,chen2019suboptimality}.

It is important that this two-replica model satisfies the coordinatewise convexity hypothesis used in the vector-spin interpolation bound.
For every $p$, the coefficient vector of the two-replica Hamiltonian is $\gamma_p(1,1)$, so its matrix-valued mixture is
\[
\boldsymbol\xi(A)=\big(\xi(A_{ij})\big)_{1\leq i,j\leq 2}.
\]
Since $\xi$ is even,
\[
\xi''(q)=\sum_{p\in 2\bbN}p(p-1)\gamma_p^2 q^{p-2}\geq 0,
\qquad q\in[-1,1].
\]
Thus every entry of $\boldsymbol\xi$ is convex in its scalar argument, as required in \cite[Theorem 3]{auffinger2022properties}.

This two-replica extension of the Parisi formula is given by a similar variational problem as in Proposition~\ref{prop:parisi-zero-temp}.
We choose suitable interpolation parameters to upper bound $\GS_{2,\eps}(\xi)$ for small $\eps$, and thus deduce non-existence of upward outlier eigenvalues at approximate ground states.
The general formula replaces $L$ by a $2\times 2$ positive semi-definite matrix $\bL$, and $\zeta$ by a cumulative distribution function on a monotone path of $2\times 2$ matrices in the positive semi-definite order.
The path of matrices is encapsulated by a function $\Phi:[0,2]\to \bbR^{2\times 2}$ with $\Phi(0)=0$ and $\Phi(2)=Q$ and $\Tr(\Phi(t))=t$ for all $t$, with $\Phi(t)-\Phi(s)$ positive semi-definite for all $t\geq s$.
Meanwhile $\alpha:[0,2]\to \bbR_+$ is the associated cumulative distribution function, which can be viewed as a positive measure on the range of $\Phi$.
For us, the relevant specialization is as follows.
Define the $2\times 2$ matrices:
    \[
    J_+
    =
    \begin{pmatrix}
    1 & 1
    \\
    1 & 1
    \end{pmatrix};
    \quad \quad 
    J_-
    =
    \begin{pmatrix}
    1 & -1
    \\
    -1 & 1
    \end{pmatrix};
    \quad \quad 
    Q
    =
    \begin{pmatrix}
    1 & 1-\eps
    \\
    1-\eps & 1
    \end{pmatrix}
    =
    \Big(1-\frac{\eps}{2}\Big)J_+ 
    + 
    \frac{\eps J_-}{2}. 
    \]
We take as given \rev{the minimizing pair} $(\zeta,L)=\argmin_{(\zeta,L)\in\cK}\cQ(\zeta,L;\xi)$ \rev{for the original model}.
Setting for convenience $\ell=\xi''(1)^{-1/2}$, their two-replica generalizations are:
    \begin{equation}
    \label{eq:two-replica-order-parameters}
    \bL
    =
    \frac{L J_+ + \ell J_-}{2}
    ;
    \quad\quad
    \alpha(t)
    =
    \zeta(t/2)/2;
    \quad\quad 
    \Phi(t)
    =
    \begin{cases}
    \frac{tJ_+}{2},\quad\quad\quad\quad\quad\quad\quad\quad~~ t\in [0,2-\eps];
    \\
    \frac{(2-\eps)J_+}{2}
    +
    \frac{\rev{(t-(2-\eps))J_-}}{2},\quad t\in [2-\eps,2].
    \end{cases}
    \end{equation}
\rev{Note that $\Phi'(t)=J_+/2$ on $t\in (0,2-\eps)$ and $\Phi'(t)=J_-/2$ on $t\in (2-\eps,2)$, while $\Phi(2)=Q$ as required.}

Applying \cite[Theorem 3]{auffinger2022properties} with $(\bL,\alpha,\Phi)$ as above gives the following upper bound on $\GS_{2,\eps}$ (see also \cite{ko2020free}).
Then Corollary~\ref{cor:constrained-GS-bound} gives a first order expansion of this bound for small $\eps$, which turns out to be sharp.

\begin{proposition}
\label{prop:2d-interpolation-bound}
    Let $\xi$ be even, with minimizer $(\zeta,L)=\argmin_{(\zeta,L)\in\cK}\cQ(\zeta,L;\xi)$. 
    Then with $\xi,\xi',\xi''$ acting entry-wise on $2\times 2$ matrices, $\la A,B\ra = \Tr(A^{\top} B)$, and $\odot$ denoting entry-wise product:
    \begin{equation}
    \label{eq:2-d-UB}
    \begin{aligned}
    2\,\GS_{2,\eps}(\xi)
    &\leq 
    \la
    \xi'(Q),\bL
    \ra
    -
    \int_0^2 
    \Big\la 
    \xi''(\Phi(t))\odot\Phi'(t),\int_0^t \alpha(s)\Phi'(s)
    \,\de s
    \rev{\Big\ra}
    \,\de t
    \\
    &\quad\quad
    +
    \int_0^2 
    \Big\la 
    \big(\bL-\int_0^t \alpha(s)\Phi'(s)\de s\big)^{-1}
    ,
    \Phi'(t)
    \Big\ra
    \,\de t
    .
    \end{aligned}
    \end{equation}
\end{proposition}

\rev{
In evaluating the above expressions, it will help to observe the simple algebraic relations
\[
J_+^2 = 2J_+,
\quad\quad 
J_-^2 = 2J_-,
\quad\quad 
J_+ J_- 
=
J_- J_+ 
=
\begin{pmatrix}
0&0\\0&0
\end{pmatrix}.
\]
Since $J_+ + J_-=2I_2$ we obtain
\[
(AJ_+ + BJ_-)^{-1}
=
(A^{-1}J_+ + B^{-1}J_-)/4.
\]
In particular, since $\la J_+,J_+\ra=4$ this yields the useful computation:
\begin{equation}
\label{eq:2-matrix-inverse}
\Big\la (AJ_+ + BJ_-)^{-1},J_+/2
\Big\ra
= 
\Big\la (A^{-1}J_+ + B^{-1}J_-)/4
,
J_+/2\Big\ra
=
\frac{1}{2A}
,\quad\forall A,B\neq 0.
\end{equation}
Similarly we have $\la J_-,J_-\ra=4$ and
$\Big\la (AJ_+ + BJ_-)^{-1},J_-/2
\Big\ra=\frac{1}{2B}$.
}

\begin{corollary}
\label{cor:constrained-GS-bound}
    For any $\iota>0$ there exists $\eps_*=\eps_*(\xi,\iota)>0$ such that for all $\eps\in (0,\eps_*)$:
    \[
    \frac{\GS_{2,\eps}(\xi)
    -
    2\,\GS(\xi)}{\eps}
    -
    \iota
    \leq 
    \sqrt{\xi''(1)}
    -
    \frac{r(\xi)}{2}
    =
    -
    \frac{\hat\zeta(1)}{2}
    \lt(
    \sqrt{\xi''(1)}-\hat\zeta(1)^{-1}
    \rt)^2
    .
    \]
\end{corollary}

\begin{proof}
    We evaluate the upper bound \eqref{eq:2-d-UB} for $2\,\GS_{2,\eps}(\xi)$ in Proposition~\ref{prop:2d-interpolation-bound} to first order in $\eps$.
    The first term is 
    \[
    L\big(\xi'(1)+(1-\eps)\xi'(1)\big)
    \rev{
    +
    \ell\big(\xi'(1)-\xi'(1-\eps)\big)
    }
    =
    2L\xi'(1)
    -
    L\xi''(1)\eps
    +
    \ell\xi''(1)\eps
    +
    O(\eps^2).
    \]
    Via the linear changes of variable $(u,v)=(s/2,t/2)$, the integral contribution from $[0,2-\eps]$ in the second term is:
    \begin{align}
    \notag
    \int\limits_0^{2-\eps}
    \lt\la 
    \frac{\xi''(t/2)J_+}{2}
    ,
    J_+
    \int_0^t 
    \zeta(s/2)/4~\de s
    \rt\ra
    \de t
    &
    \rev{=
    \frac{1}{2}
    \int\limits_0^{2-\eps}
    \xi''(t/2)\int_0^t \zeta(s/2)
    ~\de s\,\de t
    }
    \\
    \label{eq:2nd-term-GS2}
    &=
    2\int\limits_0^{1-\frac{\eps}{2}}
    \xi''(\rev{v})
    \lt( \int_0^v \zeta(u) \,\de u \rt)
    \de v
    \\
    \notag
    &=
    2\int_0^{1}
    \xi''(\rev{v})
    \lt( \int_0^v \zeta(u) \,\de u\rt)
    \de v
    -
    \lt(\xi''(1)\int_0^1 \zeta(u)\,\de u \rt)
    \eps
    +
    O(\eps^2).
    \end{align}
    The integral contribution from $[2-\eps,2]$ in the second term is $O(\eps^2)$, since the integrand takes the form 
    \[
    \Big\la 
    O(\eps)J_+ + O(1)J_-
    ,~
    O(1) J_+ + O(\eps)J_-
    \Big\ra
    =
    O(\eps).
    \]
    For the third term, we have $\Phi'(t)=J_+/2$ on $t\in [0,2-\eps]$.
    Recalling \eqref{eq:2-matrix-inverse}, we thus only need to track the $J_+$ component in the matrix $\bL-\int_0^t \alpha(s)\Phi'(s)\de s$ to determine the integral contribution on $[0,2-\eps]$.
    This yields:
    \begin{align*}
    \rev{
    \frac{1}{2}
    \int_0^{2-\eps}
    \lt(
    \frac{L}{2}-
    \frac{1}{2}
    \int_0^{t}
    \alpha(s)/2
    ~\de s
    \rt)^{-1}
    ~
    \de t
    }
    &=
    \int_0^{2-\eps}
    \lt(
    L
    -
    \int_0^{t/2}
    \zeta(u)\de u
    \rt)^{-1}
    ~
    \de t
    \\
    &=
    2
    \int_0^1
    \frac{\de q}{
    L-\int_0^{q} \zeta(u)\de u
    }
    -
    \frac{\eps}{\hat\zeta(1)}
    +
    O(\eps^2)
    .
    \end{align*}
    \rev{Finally for the contribution on $[2-\eps,2]$ in the third term, we similarly track the $J_-$ component in $\bL-\int_0^t \alpha(s)\Phi'(s)\de s$:
    \[
    \frac{1}{2}
    \int_{2-\eps}^2
    \Big(
    \frac{\ell}{2}
    -
    O(\eps)
    \Big)^{-1}
    \de t
    =
    \frac{\eps}{\ell}+O(\eps^2)
    .
    \]
    }
    Recall that $\ell=\xi''(1)^{-1/2}$, \rev{so $\ell \xi''(1)\eps+\frac{\eps}{\ell}=2\xi''(1)^{1/2}\eps$}.
    Combining the terms above and recalling that the \rev{second term contribution from \eqref{eq:2nd-term-GS2} is negated in \eqref{eq:2-d-UB}, we obtain as desired:} 
    \[
    2\,\GS_{2,\eps}(\xi)
    \leq 
    4\,\GS(\xi)
    -
    \hat\zeta(1)
    \lt(
    \sqrt{\xi''(1)}-\hat\zeta(1)^{-1}
    \rt)^2
    \eps 
    +
    O(\eps^2).
    \qedhere
    \]
\end{proof}

Using Corollary~\ref{cor:constrained-GS-bound}, we now deduce the ``no upward outliers'' property \eqref{eq:top-edge-formula} for even models in \rev{Theorem~\ref{thm:no-upward-outliers-even-xi}}.
The idea is that any upward outlier eigenvalue yields a counterexample to the estimate we have just established.

\begin{proof}[Proof of \rev{Theorem~\ref{thm:no-upward-outliers-even-xi}}, Eq.~\eqref{eq:top-edge-formula}]
    We restrict to the event of Proposition~\ref{prop:gradients-bounded}.
    Take $\eps$ small depending on $(\xi,\iota)$, and $\delta\leq \eps^2$.
    We consider the event that
    \begin{equation}
    \label{eq:event-GS-typical-v2}
    |\rev{\GS_N}-\GS(\xi)|\leq \delta
    \quad\quad
    \text{and}
    \quad\quad 
    \GS_{2,\eps}(H_N)\leq 
    2\,\GS(\xi)
    -
    \frac{\hat\zeta(1)}{2}
    \lt(
    \sqrt{\xi''(1)}-\hat\zeta(1)^{-1}
    \rt)^2
    \eps 
    +
    \iota^2\eps
    .
    \end{equation}
    This event clearly has probability at least $1-e^{-cN}$, and we will show it implies the desired conclusion.
    Let $\bsig$ be a $\delta$-approximate ground state, and suppose for sake of contradiction that
    \[
    \blambda_1(\nabla_{\sph}^2(H_N(\bsig)))
    \geq 
    -\hat\zeta(1)
    \lt(
    \sqrt{\xi''(1)}-\hat\zeta(1)^{-1}
    \rt)^2
    +
    \iota.
    \]
    (The lower bound on $\blambda_1$ is trivial since $\blambda_1(\cdot)\geq \blambda_{\lfloor\delta N\rfloor}(\cdot)$, \rev{so we just focus on the upper bound}.)
    Let 
    \[
    \theta(\eps)=\arcsin(\sqrt{\eps/2})= \sqrt{\eps/2} + O(\eps).
    \]
    With $\bv/\sqrt{N}$ the maximum unit eigenvector of $\nabla_{\sph}^2(H_N(\bsig))$, set
    \begin{equation}
    \label{eq:bsig-pm}
    \bsig^{\pm}
    =
    \cos(\theta(\eps))\bsig 
    \pm
    \sin(\theta(\eps))\bv.
    \end{equation}
    Then $R(\bsig^+,\bsig^-)=\cos(2\theta(\eps))=1-\eps$.
    Since $\eps\ll\iota$, Taylor expanding (\rev{via \eqref{eq:taylor-expansion-HN}}) gives
    \begin{align*}
    \frac{H_N(\bsig^+)
    +
    H_N(\bsig^-)
    -
    2H_N(\bsig)}{N}
    &=
    \frac{\eps}{2}
    \cdot 
    \blambda_1(\nabla_{\sph}^2(H_N(\bsig)))
    +
    O(\eps^2)
    \\
    &\geq
    -
    \frac{\hat\zeta(1)}{2}
    \lt(
    \sqrt{\xi''(1)}-\hat\zeta(1)^{-1}
    \rt)^2
    \eps 
    +
    \Omega(\iota\eps)
    .
    \end{align*}
    However, \eqref{eq:event-GS-typical-v2} implies the left-hand expression is at most
    $-
    \frac{\hat\zeta(1)}{2}
    \big(
    \sqrt{\xi''(1)}-\hat\zeta(1)^{-1}
    \big)^2
    \eps 
    +
    \iota^2 \eps.
    $
    This is a contradiction, completing the proof (with $\iota$ here functioning as $\eps$ in the original statement).
\end{proof}

\rev{
\begin{remark}
Our motivation for the order parameters in \eqref{eq:two-replica-order-parameters} was that an optimal choice $(\bsig,\bsig')$ in \eqref{eq:GS-2-eps} should essentially first decide the midpoint $\frac{\bsig+\bsig'}{2}$, and then choose the increment $\bsig-\bsig'$ to be a top eigenvector of the Hessian evaluated there.
For example when near-ground states are \emph{not} marginally stable, one might expect a positive probability for two low-temperature Gibbs samples for the pair-Hamiltonian \eqref{eq:GS-2-eps} to agree almost exactly in the first stage, but have orthogonal increments in the second, thus yielding a multiple of $J_+$ for the overlap matrix. 
This suggests that the overlap distribution should be supported on multiples of $J_+$ until the second stage.
The construction \eqref{eq:two-replica-order-parameters} accordingly behaves as a concatenation of these two stages. 
However while this suggests the general form of the construction, precisely identifying the correct order parameters required some trial and error.
\end{remark}
}

\subsection{Proof of Corollary~\ref{cor:disorder-chaos}}
\label{subsec:chaos}

Here we explain how to deduce transport disorder chaos from Corollary~\ref{cor:deep-level-sets}.

\begin{proof}[Proof of Corollary~\ref{cor:disorder-chaos}]
The proof essentially follows from \cite[Theorem 5.1]{alaoui2023shattering}, which deduces transport disorder chaos from \emph{shattering}.
Let us summarize the argument.
For $\delta\leq\delta_*(\xi)$, Corollary~\ref{cor:deep-level-sets} provides a decomposition $\{\cC_{-K},\dots,\cC_{-1},\cC_1,\dots,\cC_K\}$ of $A_{\delta}$ into separated clusters, where without loss of generality $\cC_{-i}=-\cC_i$ are antipodal pairs since $\xi$ is even.
\rev{(So far we have not asserted anything about the value of $K$.)}
The clusters $\cC_i$ have small diameter $O(\sqrt{N\delta})$ and are pair-wise separated by a much larger distance $\sqrt{N}/C(\xi)$.
Further, Lemma~\ref{lem:small-values-low-prob} implies that
\begin{equation}
\label{eq:clustering}
\mu_{\beta}\big(\bigcup_i \cC_i\big)\geq 1-e^{-cN}
\end{equation}
when $\beta\geq\beta(\xi,\delta)$ is very large.

The idea of \cite[Theorem 5.1]{alaoui2023shattering} is that $\mu_{\beta,\eps_N}$ will preserve the clustering (i.e. the property \eqref{eq:clustering}), but injects noise into the cluster \emph{weights}. 
This is implemented in Proposition 5.4 therein, by proving anti-concentration of the log-ratios $\log\big(\mu_{\beta,\eps_N}(\cC_i)/\mu_{\beta,\eps_N}(\cC_j)\big)$ for each distinct pair $i\neq j$.
The anti-concentration comes from the contribution of $\bv^{\otimes p}$ in $\tilde H_N$, where $\gamma_p>0$ and $\bv$ is chosen such that $\pm\cC_i,\pm\cC_j$ remain separated when projected along $\bv$.
The only input hypothesis that is unavailable in our setting is denoted $\mathsf{S1}$ therein: it is \textbf{not} true that $\mu_{\beta}(\cC_i)$ is exponentially small for each $i$. 
However, the proof (in the case of even models) applies with only cosmetic changes as long as $\bbE[\max_i \mu_{\beta}(\cC_i\cup \cC_{-i})]\leq 1-C(\xi)^{-1}$, i.e.\ with uniformly positive probability, the maximum probability of any cluster pair is bounded away from $1$. 
For this, it is sufficient to show that, if $\bsig,\bsig'\sim\mu_{\beta}$ are i.i.d.\ Gibbs samples, then 
\begin{equation}
\label{eq:multiple-clusters-even-generic}
\liminf_{N\to\infty}~
\bbE^{H_N}
\bbP^{\bsig,\bsig'}[|R(\bsig,\bsig')|\leq 1/2]
>0.
\end{equation}
Indeed, $|R(\bsig,\bsig')|\leq 1/2$ implies $\bsig,\bsig'$ are not in the same cluster, nor in antipodally opposite clusters $\cC_i$ and $\cC_{-i}$.

To show \eqref{eq:multiple-clusters-even-generic} we use the ``even generic'' hypothesis $\sum_{p\in 2\bbN} \frac{1_{\gamma_p\neq 0}}{p}=\infty$. 
It is standard from \cite[Chapter 3.7]{panchenko2013sherrington}, see also \cite[Section 4]{panchenko2016chaos}, that under this hypothesis, the law of $|R(\bsig,\bsig')|$ (averaged over the disorder $H_N$) has an $N\to\infty$ limit $\zeta_{\beta}$ given by the minimizer in the Parisi formula at inverse temperature $\beta$.
The precise definition of $\zeta_{\beta}$ is recalled in Section~\ref{sec:alternate-proof}, but the only property we need is the well-known fact that $0\in\supp(\zeta_{\beta})$ whenever $\xi'(0)=0$ \rev{(recall that $\supp(\cdot)$ was defined in Proposition~\ref{prop:cs-extremality-0temp})}.
In particular, \eqref{eq:multiple-clusters-even-generic} holds, so the proof of \cite[Theorem 5.1]{alaoui2023shattering} applies and yields Corollary~\ref{cor:disorder-chaos}.
\end{proof}

\section{Alternate Proof of Marginal Stability for Generic Models}
\label{sec:alternate-proof}

Here we give a different proof that full RSB near overlap $1$ implies marginal stability at low temperature.
This means we consider the typical behavior of $\bsig_{\beta}$ drawn from the Gibbs measure $\mu_{\beta}=\mu_{\beta,H_N}$ defined by
\[
\de\mu_{\beta}(\bsig)=e^{\beta H_N(\bsig)}\de\mu_0(\bsig)/Z_{N,\beta}.
\]
Here $Z_{N,\beta}=\int_{\cS_N} e^{\beta H_N(\bsig)} \de \mu_0(\bsig)$ is the partition function relative to uniform measure $\mu_0$.
The argument is based on ultrametricity of Gibbs measures, which requires $\xi$ to be \emph{generic}, i.e.
\[
    \sum_{\text{odd}~p} \frac{1_{\gamma_p\neq 0}}{p}
    =
    \sum_{\text{even}~p} \frac{1_{\gamma_p\neq 0}}{p}
    =
    \infty.
\]
(Note the similarity to the ``even generic'' hypothesis of Corollary~\ref{cor:disorder-chaos}, which would also suffice here.)
We prove the following, which is a special case of \eqref{eq:radial-deriv-formula} and \eqref{eq:marginal-stability}.\footnote{
\rev{The hypotheses are slightly stronger than our earlier full RSB endpoint assumption by Corollary~\ref{cor:full-RSB-justification}, and $\bsig_{\beta}\sim\mu_{\beta}$ is a near-ground state with high probability for large $\beta$ by Lemma~\ref{lem:small-values-low-prob}.}
}

\begin{proposition}
    \label{prop:main-generic}
    Suppose $\xi$ is generic, $1$ is not an isolated point in $T$ (recall \eqref{eq:T-def}), and $q\mapsto\xi''(q)^{-1/2}$ is non-affine on $(0,1]$. Then for any $\eps>0$, if $\bsig_{\beta}\sim\mu_{\beta}$ for $\beta\geq \beta_*(\xi,\eps)$ sufficiently large, and $\delta,c$ are small depending on $(\xi,\eps,\beta)$, and $N$ is sufficiently large, we have with probability \rev{ at least $1-\eps$}:
    \begin{align}
    \label{eq:radial-deriv-formula-generic}
        \|\nabla H_N(\bsig_{\beta})-2\sqrt{\xi''(1)}\,\bsig_{\beta}\|
        &\leq \eps\sqrt{N}
        ,
        \\
    \label{eq:bulk-edge-formula-generic}
        |\blambda_1\big(\nabla_{\sph}^2 H_N(\bsig_{\beta})\big)|
        +
        |\blambda_{\lfloor\delta N\rfloor}\big(\nabla_{\sph}^2 H_N(\bsig_{\beta})\big)|
        &\leq \eps
        .
    \end{align}
\end{proposition}

The positive-temperature Parisi formula will enter only by prescribing the possible overlaps of Gibbs samples, giving a more conceptual explanation for this implication.
We now recall its statement, which gives the limiting value of the free energy $F_{N,\beta}=\frac{1}{N} \log Z_{N,\beta}$.

Let $\cM$ denote the space of all right-continuous non-decreasing functions $\rev{\zeta_{\beta}}:[0,1]\to [0,1]$ with $\rev{\zeta_{\beta}}(\hat q)=1$ for some $\hat q<1$ (which may depend on $\rev{\zeta_{\beta}}$).
\rev{(Here we have in mind the positive temperature setting $\beta<\infty$, but the definition of the space $\cM$ does not depend on $\beta$.)}
Let $\rev{\hat\zeta_{\beta}}(q) = \int_q^1 \rev{\zeta_{\beta}}(s)~\de s$ and define the Crisanti--Sommers functional
\begin{equation}
    \label{eq:cs-functional}
    \cP(\rev{\zeta_{\beta}};\xi)
    = \fr12 \Big(
        \xi'(0) \rev{\hat\zeta_{\beta}}(0)
        + \int_0^1 \xi''(q)\rev{\hat\zeta_{\beta}}(q) ~\de q
        + \int_0^\hq \fr{\de q}{\rev{\hat\zeta_{\beta}}(q)}
        + \log (1-\hq)
    \Big).
\end{equation}
Note that $\rev{\hat\zeta_{\beta}}(q) = 1-q$ for $q > \hq$, so this functional is independent of $\hq$.
The spherical Parisi formula at positive temperature is as follows.

\begin{proposition}[{\cite{talagrand2006free,chen2013aizenman,chen2017parisi}}]
\label{prop:parisi}
    For $\beta\in \bbR_+$, the asymptotic free energy satisfies:
    \begin{align*}
    F(\beta)
    &\equiv
    \plim_{N\to\infty} F_{N,\beta}
    =
    \rev{
    \inf_{\zeta_{\beta}\in\cM} \cP(\zeta_{\beta};\beta\xi)},
    \\
    \GS(\xi)
    &=
    \lim_{\beta\to\infty} F(\beta)/\beta
    .
    \end{align*}
    Further \rev{a unique minimizer exists, which we also denote by $\zeta_{\beta}$}. With $(\zeta,L)$ the zero-temperature minimizers, one has
    \[
    \zeta =
    \lim_{\beta\to\infty} \beta \zeta_{\beta},
    \quad\quad
    L=\lim_{\beta\to\infty} \int_0^1 \beta \zeta_{\beta}(s)~\de s.
    \]
    Here $\zeta$ and $\zeta_{\beta}$ are metrized by the vague topology on $[0,1)$ for the corresponding positive measures as in \eqref{eq:def-nu-infty}.
\end{proposition}

There is a positive-temperature analog of Proposition~\ref{prop:cs-extremality-0temp} characterizing $\rev{\zeta_{\beta}}$, but this will not be needed.
We instead use the following qualitative result on the overlaps between Gibbs samples.
Recall that $\supp(\zeta_{\beta})\subseteq [0,1)$ is the set of points of increase of $\zeta_{\beta}$.

\begin{corollary}
\label{cor:full-rsb-positive-temp}
    Assume the conditions of Proposition~\ref{prop:main-generic}.
    Then there exists $\eps=\eps_*(\xi)>0$ such that for any $q\in (1-\eps_*,1)$ and $\delta>0$, if $\beta\geq \beta_*(\xi,q,\delta)$ is sufficiently large, then $\supp(\zeta_{\beta})\cap (q-\delta,q+\delta)\neq\emptyset$.
\end{corollary}

\begin{proof}
    As stated in \eqref{prop:parisi}, the zero-temperature optimizer $\zeta$ is the vague limit as $\beta\to\infty$ of $\beta\zeta_{\beta}$. 
    Hence it suffices to show the first claim that $\zeta$ is strictly increasing in a neighborhood of $1$, which follows from Corollary~\ref{cor:full-RSB-justification} and the non-affineness assumption in Proposition~\ref{prop:main-generic}.
\end{proof}

Next we define the event that a large constant number $K$ of Gibbs samples have all overlaps approximately $q$:
\begin{equation}
\label{eq:ultrametric-star-event}
    E_{K,q,\delta}\equiv 
    \lt\{
    |R(\bsig_{\beta}^i,\bsig_{\beta}^j)-q|\leq \delta~~\forall\, 1\leq i<j\leq K\rt\}.
\end{equation}
Here $\bsig_{\beta}^1,\dots,\bsig_{\beta}^K\stackrel{\iid}{\sim}\mu_{\beta,H_N}$ are always \iid samples from $\mu_{\beta}$. 
Our alternate proof of Proposition~\ref{prop:main-generic} relies on the fact that $E_{K,q,\delta}$ has uniformly positive probability when $q\in\supp(\zeta_{\beta})$, even conditional on typical $(H_N,\bsig_{\beta}^1)$.
This property follows from the Ghirlanda-Guerra identities, which hold when $\xi$ is generic.
Namely as explained in \cite[Chapter 3.7]{panchenko2013sherrington}, when $\xi$ is generic the limiting overlap arrays of i.i.d.\ Gibbs samples converge to a limiting random overlap structure as $N\to\infty$. 
It is easy to see that this limiting structure satisfies the next proposition.

\begin{proposition}
\label{prop:generic-overlaps}
    If $\xi$ is generic, then for all $q\in\supp(\zeta_{\beta})$ and any $\delta>0$,
    \[
    \lim_{\eta\to 0}
    \lim_{N\to\infty}
    \bbP
    \lt[
    \bbP[E_{K,q,\delta}~|~(H_N,\bsig_{\beta}^1)]>\eta
    \rt]
    =1.
    \]
\end{proposition}


The next proposition shows that low dimensional spaces do not contain a large number of points with all equal distances.
Though elementary, it is key to our argument.

\begin{proposition}
\label{prop:non-embed}
    For any $d>0$ there exists $\eps>0$ such that no $d+2$ points $x_1,\dots,x_{d+2}\in \bbR^d$ satisfy 
    \[
    \|x_i-x_j\|\in [a(1-\eps),a(1+\eps)],\quad\forall~1\leq i<j\leq d+2
    \]
    for any $a>0$.
\end{proposition}

\begin{proof}
    Without loss of generality set $a=1$.
    In the proof below, implicit constants in $O(\cdot)$ may depend on $d$.
    Suppose such points exist and let $z=\frac{1}{d+2}\sum_{i=1}^{d+2} x_i$ and $y_i=x_i-z$.
    Note that for each $j,k\neq i$:
    \[
    \langle x_i-x_j,x_i-x_k\rangle 
    =
    \frac{\|x_i-x_j\|^2 + \|x_i-x_k\|^2 - \|x_j-x_k\|^2}{2}
    =
    \frac{1}{2}\rev{+}O(\eps).
    \]
    We compute
    \[
    \langle y_i,y_i\rangle
    =
    (d+2)^{-2}
    \sum_{k,\ell=1}^{d+2}
    \langle x_i-x_k,x_i-x_{\ell}\rangle 
    \\
    =
    \lt(\frac{d+1}{d+2}\rt)^2 \rev{+} O(\eps).
    \]
    Similarly $|\langle x_i-x_k,x_j-x_{\ell}\rangle|\leq O(\eps)$ when all four indices are distinct, so for $i\neq j$:
    \[
    \langle y_i,y_j\rangle
    =
    (d+2)^{-2}
    \sum_{k,\ell=1}^{d+2}
    \langle x_i-x_k,x_j-x_{\ell}\rangle 
    =
    -\frac{d+1}{(d+2)^2} \rev{+} O(\eps).
    \]
    Hence the $(d+2)\times (d+2)$ matrix $M$ with entries $M_{i,j}=\langle y_i,y_j\rangle$ is entrywise within $O(\eps)$ of 
    \[
    \wt M_{i,j}= 
    \begin{cases}
    \lt(\frac{d+1}{d+2}\rt)^2,\quad i=j,
    \\
     -\frac{d+1}{(d+2)^2},\quad i\neq j.
    \end{cases}
    \]
    Diagonal dominance implies $\rank(\wt M)=d+1$, hence $\rank(M)\geq d+1$ for $\eps$ sufficiently small.
    However by construction $\rank(M)\leq d$ since $y_1,\dots,y_{d+2}\in\bbR^d$. This is a contradiction and completes the proof.
\end{proof}

\subsection{Proof of Proposition~\ref{prop:main-generic}}

In light of Lemma~\ref{lem:spectrum-approx}, it will suffice to prove \eqref{eq:bulk-edge-formula-generic} with $\lfloor \delta N\rfloor$ replaced by a large constant $K$. 
Thus, for sake of contradiction we fix $K,C>0$ such that 
\begin{equation}
\label{eq:K-c-contradiction}
    \limsup_{\beta\to\infty}
    \limsup_{N\to\infty}
    \bbP[
    \blambda_K(\nabla^2_{\sph} H_N(\bsig_{\beta}))\leq -C
    ]>0.
\end{equation}
We next choose several more constants: $\eps\leq \eps_*(\xi,K,C)$ is taken sufficiently small such that $q=1-\eps\in\supp(\zeta)$, and we set $\lambda=\eps^{0.55}$. 
Then we send $\beta\to\infty$, inducing a choice of $\eta\to 0$ so that Corollary~\ref{cor:near-ground-state} holds, and $\delta\to 0$ so Corollary~\ref{cor:full-rsb-positive-temp} holds.
Thus $\beta$ (resp. $\eta,\delta)$ is sufficiently large (resp. small) depending on $(\xi,K,C,\eps)$.

Given $\bsig\in\cS_N$, let $\cT(\bsig)$ be the tangent space to $\cS_N$, viewed as a codimension $1$ linear subspace of $\bbR^N$.
Let $U_K=U_K(\bsig)\subseteq\cT(\bsig)$ denote the span of the top $K$ eigenvectors of $\nabla_{\sph}^2 H_N(\bsig)$, and let $U_K^{\perp}(\bsig)\subseteq\cT(\bsig)$ be its orthogonal complement in $\cT(\bsig)$.
For any $\bv\in\bbR^N$, there is a unique decomposition
\begin{equation}
    \label{eq:parallel-perp-decomp}
    \bv=\bv_{\parallel}+\bv_{\perp}+R(\bsig,\bv)\bsig
\end{equation}
with $(\bv_{\parallel},\bv_{\perp})\in U_K\times U_K^{\perp}$. Note that if $\bv=\bsig'-\bsig$ then $\bv_{\parallel}+\bv_{\perp}$ is proportional to the derivative $\gamma'(0)$ for $\gamma$ a geodesic path from $\bsig$ to $\bsig'$.
We also set $U_{K,\lambda}$ to be the $\lambda\sqrt{N}$-neighborhood of $U_K(\bsig)+\bsig$ in $\bbR^N$ (note that $U_K\subseteq\cT(\bsig)\subseteq\bbR^N$ can be naturally viewed as a \rev{linear subspace} of $\bbR^N$).

\begin{proposition}
\label{prop:taylor-expansion}
    For $\bsig,\bsig'\in\cS_N$, let $\bv=\bsig'-\bsig$ and define the decomposition \eqref{eq:parallel-perp-decomp} based on $U_K(\bsig)$. Then uniformly over $H_N\in K_N$:
    \[
    H_N(\bsig')= H_N(\bsig) 
    + 
    \la \nabla_{\sph} H_N(\bsig),\bv_{\parallel}+\bv_{\perp}\ra 
    +
    \la \nabla_{\sph}^2 H_N(\bsig),(\bv_{\parallel}+\bv_{\perp})^{\otimes 2}\ra
    +
    N^{-1/2}\cdot O(\|\bsig'-\bsig\|^3)
    .
    \]
\end{proposition}

\begin{proof}
    \rev{Recalling \eqref{eq:taylor-expansion-HN}}, this amounts to a second order Taylor expansion of $H_N$ along the geodesic path from $\bsig$ to $\bsig'$. The required estimate on the third derivative holds since $H_N\in K_N$.
\end{proof}

We now show the Gibbs mass near $\bsig_{\beta}$ essentially lives within the set $U_{K,\lambda}(\bsig_{\beta})$.

\begin{lemma}
\label{lem:local-subspace-concentration}
    For constants chosen as above, suppose $H_N\in K_N$ and that $\bsig\in\cS_N$ satisfies:
    \begin{align*}
        \|\nabla H_N(\bsig)\|
        \leq 
        \eta\sqrt{N};
        \quad\quad
        \blambda_1(\nabla^2_{\sph} H_N(\bsig))\leq \eta;
        \quad\quad
        \blambda_K(\nabla^2_{\sph} H_N(\bsig))\leq -C.
    \end{align*}
    Then it follows that
    \[
    \mu_{\beta}
    \big(B_{2\sqrt{\eps N}}(\bsig)\backslash U_{K,\lambda}(\bsig)
    \big)
    \leq 
    e^{-cN}
    .
    \]
\end{lemma}

\begin{proof}
    Given $\bsig'\in B_{2\sqrt{\eps N}}(\bsig)$, let $\bv=\bsig'-\bsig$ and let $\bw=\bv_{\parallel}+\bv_{\perp}$ be the associated tangent vector to $\cS_N$ at $\bsig$ with $\|\bw\|\leq 2\sqrt{\eps N}$.
    From the way we chose constants, $\lambda^2$ is larger than $\eta\eps^{1/2}+\eps^{3/2}$ by a super-constant factor. 
    Hence Proposition~\ref{prop:taylor-expansion} yields
    \begin{align*}
    H_N(\bsig')
    &\leq 
    H_N(\bsig)
    +
    \eta 
    (\|\bw\|\sqrt{N}+\|\bw\|^2)
    -
    C\|\bv^{\perp}\|^2/2
     +
    O(\eps^{3/2} N)
    \\
    &\leq
    H_N(\bsig) + O(\eta \eps^{1/2}+\eps^{3/2}) N
    -
    C\|\bv^{\perp}\|^2/2
    \\
    &
    \leq
    H_N(\bsig)-
    C\lambda^2 N/3.
    \end{align*}
    Lemma~\ref{lem:small-values-low-prob} now completes the proof.
\end{proof}

\begin{lemma}
\label{lem:embedding-failure}
    With parameters as above (in particular $\delta$ small depending on $\eps$), suppose $\bsig^1,\dots,\bsig^{K+3}\in\cS_N$ satisfy:
    \begin{align*}
        R(\bsig^i,\bsig^j)
        &\in [1-\eps-\delta,1-\eps+\delta],
        \quad \forall ~1\leq i<j\leq K+3
        .
    \end{align*}
    Then with decomposition~\eqref{eq:parallel-perp-decomp} defined based on $U_K(\bsig^1)$,
    and with $\bv^j=\bsig^j-\bsig^1$ for $2\leq j\leq K+3$,
    \[
    \max_{2\leq j\leq K+3} \|\bv^j_{\perp}\|> \lambda\sqrt{N}.
    \]
\end{lemma}

\begin{proof}
    Suppose not. Then the vectors $\bv^j_{\parallel}$ for $2\leq j\leq K+3$
    lie in a $K$-dimensional subspace and for $i\neq j$:
    \begin{align*}
    \|\bv^i_{\parallel}-\bv^j_{\parallel}\|
    &=
    \|\bv^i-\bv^j\|
    \rev{+} 
    O(\|\bv^i-\bv^i_{\parallel}\|+\|\bv^j-\bv^j_{\parallel}\|)
    \\
    &=
    \|\bsig^i-\bsig^j\|
    \rev{+} 
    O(\|\bv^i-\bv^i_{\parallel}\|+\|\bv^j-\bv^j_{\parallel}\|)
    \\
    &=
    \Big(\sqrt{2-2(1-\eps)^2} \rev{+} O(\delta \eps^{-1/2} + \lambda)\Big)\sqrt{N}
    .
    \end{align*}
    This is impossible by Proposition~\ref{prop:non-embed} since $\sqrt{2-(2-\eps)^2}\asymp\sqrt{\eps}$
    and $\delta\eps^{-1/2}+\lambda\leq o(\sqrt{\eps})$, completing the proof.
\end{proof}

\begin{proof}
[Proof of Proposition~\ref{prop:main-generic}]
    First, Corollary~\ref{cor:near-ground-state} implies $\|\nabla_{\sph}H_N(\bsig_{\beta})\|\leq \eps\sqrt{N}/2$.
    It suffices to establish \eqref{eq:bulk-edge-formula-generic} with $\lfloor \delta N\rfloor$ replaced by $K$, as Lemma~\ref{lem:spectrum-approx} shows these are equivalent and also then yields \eqref{eq:radial-deriv-formula-generic}.
    Combining Lemma~\ref{lem:spectrum-approx} with \eqref{eq:near-concave} gives the upper bound for $\blambda_1\big(\nabla_{\sph}^2 H_N(\bsig_{\beta})\big)$. 
    
    The main part of the proof is the lower bound on $\blambda_K(\nabla^2_{\sph} H_N(\bsig_{\beta}))$. 
    Suppose for sake of contradiction that \eqref{eq:K-c-contradiction} holds for some $K,C>0$, and define the event 
    \[
    S_{\eigen}
    =
    \Big\{\blambda_K(\nabla^2_{\sph} H_N(\bsig_{\beta}))\leq -C\Big\}. 
    \]
    Choose small $\eps$ depending on $(K,C)$ such that $1-\eps\in\supp(\zeta)$. 
    Let $\eta\ll \delta$ be small depending on $(\xi,K,C,\eps,\beta)$, and define for \iid Gibbs samples $(\bsig_{\beta}^1,\dots,\bsig_{\beta}^{K+3})$ the event
    \[
    S_{\generic}
    =
    \Big\{\bbP[E_{K+3,1-\eps,\delta}~|~(H_N,\bsig_{\beta}^1)]>(K+2)\eta\Big\}.
    \]
    By Proposition~\ref{prop:generic-overlaps}, for 
    $N\geq N_0(\xi,K,C,\eps,\beta,\delta,\eta)$ sufficiently large we have
    \[
    \bbP
    \lt[
    S_{\generic}
    \rt]
    >1-\frac{\eps}{2}.
    \]
    The event $E_{K+3,1-\eps,\delta}$ trivially implies $\bsig_{\beta}^i\in B_{2\sqrt{\eps N}}(\bsig_{\beta}^1)$ for all $2\leq i\leq K+3$ because $\delta$ is small depending on $\eps$.
    Explicitly, if $R(\bsig,\bsig')\geq 1-\eps-\delta$ then $\|\bsig-\bsig'\|=\sqrt{2(\eps+\delta) N}\leq 2\sqrt{\eps N}$.
    Let $S_{\bdd}$ be the event $H_N\in K_N$, which has probability $1-e^{-cN}$ by Proposition~\ref{prop:gradients-bounded}.
    We claim the three events $S_*$ cannot all hold, i.e. deterministically,
    \begin{equation}
    \label{eq:S-intersection-empty}
        S_{\generic}\cap S_{\bdd}\cap S_{\eigen}=\emptyset.
    \end{equation}
    Indeed, assume that $S_{\generic}$ and $S_{\bdd}$ hold. 
    Then using Lemma~\ref{lem:embedding-failure} in the first step and the definition of $S_{\generic}$,
    \[
    \mu_{\beta}\big(B_{2\sqrt{\eps N}}(\bsig_{\beta}^1)\backslash U_{K,\lambda}(\bsig^1)\big) 
    \geq 
    \bbP[E_{K+3,1-\eps,\delta}~|~(H_N,\bsig_{\beta}^1)]/(K+2)\geq \eta.
    \]
    In light of Lemma~\ref{lem:local-subspace-concentration} and the assumption $S_{\bdd}$, we find that $S_{\eigen}$ indeed cannot hold. 
    Finally \rev{having shown \eqref{eq:S-intersection-empty}}, it follows that $\bbP[S_{\eigen}]\leq \frac{\eps}{2}+e^{-cN}\leq \eps$ which concludes the proof.
\end{proof}

\subsection*{Acknowledgements}

Thanks to Brice Huang for helpful discussions, and for comments on a previous version of the paper.
We sincerely thank Qiang Zeng for pointing out the affine exception in Corollary~\ref{cor:full-RSB-justification} and that an earlier argument for the absence of \emph{downward} outliers was incorrect due to its reliance on interpolation bounds for a nonconvex vector spin glass.
We also thank Andrea Montanari, Pierfrancesco Urbani, and Lenka Zdeborov{\'a} for enlightening conversations, as well as the anonymous referee for numerous helpful suggestions.



\small

\bibliographystyle{alphaabbr}
\bibliography{bib}

\end{document}